\documentclass[12pt]{amsart}
\usepackage{latexsym}
\usepackage{amssymb}
\usepackage{amscd}
\usepackage{mathrsfs}
\renewcommand\a{\alpha}

\renewcommand\d{\delta}
\newcommand\la{\lambda}

\renewcommand\th{\theta}

\newcommand\io{\iota}

\newcommand\s{\sigma}

\renewcommand\t{\tau}

\newcommand\vL{\varLambda}

\newcommand\vG{\varGamma}
\newcommand\ve{\varepsilon}

\newcommand{\OO}{\mathbb O}

\newcommand\SB{\mathscr{B}}
\newcommand\SC{\mathscr{C}}

\newcommand\SE{\mathscr{E}}

\newcommand\SN{\mathscr{N}} 
\newcommand\SO{\mathscr{O}}
\newcommand\SP{\mathscr{P}}
\newcommand\SQ{\mathscr{Q}}

\newcommand\SH{\mathscr{H}}
\newcommand\ScS{\mathscr{S}}
\newcommand\SW{\mathscr{W}}

\newcommand\SX{\mathscr{X}}

\newcommand\Ql{\bar{\mathbf Q}_l}

\newcommand\BP{\mathbf P}
\newcommand\BQ{\mathbf Q}

\newcommand\BC{\mathbf C}

\newcommand\BZ{\mathbf Z}

\newcommand\BV{\mathbf V}

\newcommand\Bm{\mathbf m}

\newcommand\Bv{\mathbf v}
\newcommand\Bk{\mathbf k}

\newcommand\Bla{\boldsymbol\lambda}

\newcommand\Bmu{\boldsymbol\mu}

\newcommand\Fg{\mathfrak g}

\newcommand\Fh{\mathfrak h}
\newcommand\Fs{\mathfrak s}
\newcommand\Fl{\mathfrak l}

\newcommand\Fp{\mathfrak p}

\newcommand\iv{^{-1}}
\newcommand\wh{\widehat}
\newcommand\wt{\widetilde}
\newcommand\wg{^{\wedge}}

\newcommand\ol{\overline}

\newcommand\hra{\hookrightarrow}
\newcommand\lra{\leftrightarrow}

\newcommand\IC{\operatorname{IC}}

\newcommand\End{\operatorname{End}}

\newcommand\ind{\operatorname{ind}}

\newcommand\Lie{\operatorname{Lie}}

\newcommand\ch{\operatorname{ch}}

\newcommand\uni{_{\operatorname{uni}}}
\newcommand\nil{_{\operatorname{nil}}}

\newcommand\lp{\operatorname{\!\langle\!}}
\newcommand\rp{\operatorname{\!\rangle\!}}

\newcommand\even{\operatorname{even}}
\newcommand\odd{\operatorname{odd}}

\newcommand{\isom}{\,\raise2pt\hbox{$\underrightarrow{\sim}$}\,}
\numberwithin{equation}{section}

\newtheorem{thm}{Theorem}[section]
\newtheorem{lem}[thm]{Lemma}
\newtheorem{cor}[thm]{Corollary}
\newtheorem{prop}[thm]{Proposition}

\def \para#1{\par\medskip\textbf{#1}
              \addtocounter{thm}{1}}

\def \remark#1{\par\medskip\noindent
                \textbf{Remark #1}
                \addtocounter{thm}{1}}

\def \remarks#1{\par\medskip\noindent
                \textbf{Remarks #1}
                \addtocounter{thm}{1}}

\begin{document}
\setlength{\baselineskip}{4.9mm}
\setlength{\abovedisplayskip}{4.5mm}
\setlength{\belowdisplayskip}{4.5mm}
\renewcommand{\theenumi}{\roman{enumi}}
\renewcommand{\labelenumi}{(\theenumi)}
\renewcommand{\thefootnote}{\fnsymbol{footnote}}
\renewcommand{\thefootnote}{\fnsymbol{footnote}}
\allowdisplaybreaks[2]
\parindent=20pt
\medskip
\begin{center}
 {\bf Springer correspondence for symmetric spaces }
\\
\vspace{1cm}
Toshiaki Shoji
\\ 
\title{}
\end{center}

\begin{abstract}
This is a survey on the Springer correspondence for 
symmetric spaces.  We discuss various generalization  
of the theory of the Springer correspondence for reductive 
groups to symmetric spaces and exotic symmetric spaces 
associated to classical groups.  
\end{abstract}

\maketitle
\pagestyle{myheadings}

\par\noindent
{\bf \S 1. Introduction }
\par\medskip
The Springer correspondence is the canonical correspondence between 
the unipotent classes of a connected reductive group $G$ and 
irreducible representations of the Weyl group $W$ of $G$, 
established by Springer [Sp] in 1976.  In 1981, Lusztig [L2] 
found a way to reformulate the Springer's theory by means of the theory
of perverse sheaves.  In the same paper, he gave a geometric
interpretation of Koskta polynomials in terms of the intersection 
cohomology associated to unipotent classes of $GL_n$.
The theory of Springer correspondence was generalized by Lusztig [L3]
in 1984 to the theory of generalized Springer correspondence, which 
became a basis of his theory of character sheaves [L4]. 
\par
It is an interesting problem to extend the theory of character sheaves 
to a more general setting, such as a variety on which $G$ acts.  
Already in 1989, Ginzburg [Gi] introduced the character sheaves on 
the symmetric space $G/H$, and Grojnowski [Gr] and Henderson [H] studied 
the case where $G/H = GL_{2n}/Sp_{2n}$ extensively.    
Recently, by Achar-Henderson [AH], Finkelberg-Ginzburg-Travkin [FGT], 
Kato [K1], Shoji-Sorlin [SS], Shoji [S3],  different types of examples, 
such as the enhanced variety and the exotic symmetric spaces associated 
to symplectic groups, are found. Those examples enjoy the satisfied 
theory of the Springer correspondence, as an analogue of the Springer 
correspondence for $GL_n$. In turn, in Chen-Vilonen-Xue [CVX], Shoji-Yang [SY], 
the case of the symmetric spaces associated to orthogonal groups were studied.  
In this case, an analogue of the generalized Springer correspondence appears, 
rather than the Springer correspondence.  
We note that the theory of perverse sheaves associated to 
(Lie algebra version of) the symmetric spaces 
of general type was studied by Lusztig-Yun [LY]. 
But in this case, the Springer correspondence does not hold in general, 
in the strict sense. 
\par
This paper is a survey on the (generalized) Springer correspondence for 
symmetric spaces (mainly associated to classical groups), based on the talks 
at the conference ALTReT2019 in Ito, Japan.         
  
\par\bigskip\noindent
{\bf Contents }
\par\bigskip\noindent
\S 1. Introduction \\
\S 2. Springer correspondene for $GL_n$ \\
\S 3. Springer correspondence for reductive groups \\
\S 4. The interpretation via perverse sheaves \\
\S 5. The generalized Springer correspondence  \\
\S 6. Geometric realization of Kostka polynomials \\
\S 7. The enhanced variety $GL(V) \times V$ \\
\S 8. Springer correspondence for the enhanced variety \\
\S 9. Double Kostka polynomials \\
\S 10. Symmetric spaces in algebraic setting  \\
\S 11. Symmetric spaces associated to classical groups \\
\S 12. Unipotent orbits in $H$ and in $G^{\io\th}\uni$  \\
\S 13. Exotic symmetric space associated to symplectic groups \\
\S 14. Springer correspondence for $G^{\io\th}\uni$, the case $H = Sp_N$  \\
\S 15. Springer correspondence for $G^{\io\th}\uni$, the case $H = SO_{2n+1}$ \\
\S 16. Generalized Springer correspondence for $G^{\io\th}\uni$, 
         the case $H = SO_{2n+1}$  \\
\S 17. Exotic symmetric spaces of higher level, the case $H = Sp_N$  \\
\S 18. Exotic symmetric spaces, the case $H = SO_{2n+1}$ \\
\S 19. Symmetric spaces in characteristic 2  
  
\par\bigskip\noindent
{\bf \S 2.  Springer correspondence for $GL_n$ }
\addtocounter{section}{2}
\par\medskip
Throughout the paper, we assume that $\Bk$ is an algebraically closed field 
of characteristic $p \ge 0$, and we consider $\Ql$-sheaves, where $\Ql$ is 
an algebraic closure of the $\l$-adic number field $\BQ_l$ with $l \ne p$. 
Note that if $p = 0$, one can replace $\Bk$ by the complex number field $\BC$, 
and $\Ql$-sheaves by ordinary $\BC$-sheaves.  All the  representations 
of finite groups are considered over $\Ql \simeq \BC$.
For a finite group $\vG$, we denote by $\vG\wg$ the set of isomorphism classes 
of irreducible representations of $\vG$. 

\par  
We will start from the simplest example.  
Let $G = GL_n(\Bk)$, $B$ a Borel subgroup of $G$ containing 
a maximal torus $T$. Let $W = N_G(T)/T$ be the Weyl group of $G$. 
Hence $W$ is isomorphic to the symmetric group $S_n$ of degree $n$.   
Let $\SP_n$ be the set of partitions $\la = (\la_1, \dots, \la_k)$ 
such that $\sum_i\la_i = n$. 
It is well-known that there is a bijection 
\begin{equation*}
\tag{2.1}
S_n\wg \simeq \SP_n, \quad (V_{\la} \lra \la)
\end{equation*}
where we normalize this so that $V_{\la}$ is the identity  
(resp. the sign) representation if $\la = (n)$ (resp. $\la = (1^n)$). 
\par
On the other hand, let $G\uni$ be set of unipotent elements in $G$, 
called the {\bf unipotent variety} of $G$. 
$G$ acts on $G\uni$ by the conjugation, and $G\uni$ is a union of 
unipotent classes of $G$.  We denote by $G\uni/\!\!\sim_G$ the set of unipotent
classes in $G$.   
It is also well-known, via the Jordan normal form, that 
\begin{equation*}
\tag{2.2}
G\uni/\!\!\sim_G  \ \simeq \SP_n, \quad (\SO_{\la} \lra \la).
\end{equation*} 

It follows from (2.1) and (2.2) that there exists a bijection 
$S_n\wg \simeq G\uni/\!\!\sim_G$ by $V_{\la} \lra \SO_{\la}$. 
But this is nothing more than  the 
parametrization of two sets, $S_n\wg$ and $G\uni/\!\!\sim_G$, coincides 
each other, by accident. It would be more interesting, form a mathematical 
point of view,  to show that there exists a ``canonical'' bijection 
between $S_n\wg$ and $G\uni/\!\! \sim _G$, independent from  the parametrization. 
Actually, this assertion was achieved by the discovery of the Springer correspondence, 
as explained below. 
\par
Let $\SB = G/B$ be the flag variety of $G$. Consider the variety 

\begin{equation*}
\tag{2.3}
\wt G\uni = \{ (x, gB) \in G\uni \times \SB \mid g\iv xg \in B \}
\end{equation*}
and a map $\pi : \wt G\uni \to G\uni$, $(x, gB) \mapsto x$. 
Then $\wt G\uni$ is smooth, irreducible, and $\pi$ is proper.  In fact, 
$\pi$ gives a resolution of singularities of $G\uni$, and is called 
the {\bf Springer resolution} of $G\uni$.
For $x \in G\uni$, we define a closed set $\SB_x$ of $\SB$ by 

\begin{equation*}
\tag{2.4}
\pi\iv(x) \simeq \SB_x = \{ gB \in \SB \mid g\iv xg \in B\}.
\end{equation*}

$\SB_x$ is called the {\bf Springer fibre} of $x$.
$\SB_x$ is not smooth, nor irreducible, but it is an interesting variety. 
We consider the cohomology group $H^i(\SB_x, \Ql)$.  
Then $H^i(\SB_x, \Ql) = 0$ if $i > 2d_x$, where $d_x = \dim \SB_x$. 
Thus $H^{2d_x}(\SB_x, \Ql)$ is the top cohomology. 
The following result holds.

\begin{thm}[Springer]   
Let $x \in G\uni$. 
\begin{enumerate}
\item    $H^i(\SB_x, \Ql)$ has a structure of $S_n$-module, called 
the {\bf Springer representation} of $S_n$.
\item  $H^{2d_x}(\SB_x, \Ql)$ is an irreducible $S_n$-module. 
\item   The map $x \mapsto H^{2d_x}(\SB_x, \Ql)$ induces  
a ``canonical'' bijection 
\begin{equation*}
G\uni/\!\!\sim_G \isom S_n\wg. 
\end{equation*}
\item  For $x \in \SO_{\la}$, 
$H^{2d_x}(\SB_x,\Ql) \simeq V_{\la}$ as $S_n$-modules. 
Hence we obtain the correspondence $\SO_{\la} \lra V_{\la}$.   
\end{enumerate}
The bijective correspondence in {\rm (iii)} is called 
the {\bf Springer correspondence}.
\end{thm}
\par\bigskip\noindent
{\bf \S 3.  Springer correspondence for reductive groups }
\par\medskip
\addtocounter{section}{1}
\addtocounter{thm}{-1}

The Springer correspondence for $G  = GL_n$ gives a bijective correspondence
between unipotent classes of $G$ and irreducible representations of 
the Weyl group $W$ of $G$. But this does not hold for reductive groups in 
general, and needs some modification. 
Let $G$ be a connected reductive group. $B, T, W, G\uni$, etc. are defined 
similarly as in Section 2. The Springer resolution $\pi : \wt G\uni \to G\uni$
is defined similarly.  (Note in the case where $p > 0$, we need to assume that
$G$ is simply connected for obtaining the resolution of singularities.  But 
we ignore this point, and use the ``Springer resolution'' for reductive groups 
in general.)  The Springer fibre $\SB_x$ is defined similarly, and 
the representation of $W$ on $H^i(\SB_x, \Ql)$, namely the Springer representation 
of $W$, was constructed by Springer [Sp]. 
But the top cohomology 
$H^{2d_x}(\SB_x, \Ql)$ is not necessarily irreducible.  
\par
For $x \in G\uni$, put $A_G(x) = Z_G(x)/Z_G^0(x)$. Then $A_G(x)$ is a finite 
group.  $Z_G(x)$ acts on $\SB_x$ by $z : gB \mapsto zgB$, which induces an action 
of $Z_G(x)$ on $H^i(\SB_x, \Ql)$.  Since $Z_G^0(x)$ acts trivially on 
$H^i(\SB_x, \Ql)$, we have an action of $A_G(x)$ on $H^i(\SB_x, \Ql)$. 
It is proved that the action of $W$ on $H^i(\SB_x, \Ql)$ commutes with $A_G(x)$, 
thus we have an action of $W \times A_G(x)$ on $H^i(\SB_x, \Ql)$. 

\remark{3.1.} \ In the case where $G = GL_n$, $Z_G(x)$ is connected, hence 
$A_G(x) = \{ 1 \}$. For $G = Sp_N, SO_N$, $A_G(x) \simeq (\BZ/2\BZ)^c$ for some $c$.
In turn, for $G$ of type $G_2, F_4, E_8$, there exist a unique class $\SO$ in $G$
such that $A_G(x) \simeq S_3, S_4$ and $S_5$ for $x \in \SO$, respectively. 
\par\medskip

The following theorem  is the original form of the Springer correspondence 
due to Springer.

\begin{thm}[{Springer [Sp]}]  
Let $G$ be reductive, and $x \in G\uni$.
\begin{enumerate}
\item
Consider the decomposition of $H^{2d_x}(\SB_x, \Ql)$ 
by $A_G(x) \times W$-modules, 
\begin{equation*}
H^{2d_x}(\SB_x, \Ql) \simeq \bigoplus_{\tau \in A_G(x)\wg}\tau \otimes V_{(x,\tau)}.
\end{equation*}
Then $V_{(x,\tau)}$ is an irreducible $W$-module if it is non-zero. 
\item
All the irreducible $W$-modules are realized in this way uniquely. 
In particular, 

we have an injective map
\begin{equation*}
W\wg \hra  \SN_G := \{ (x,\tau) \mid x \in G\uni/\!\sim_G, \tau \in A_G(x)\wg \}. 
\end{equation*} 
\end{enumerate}
\end{thm}
   
\remarks{3.3.} \  
(i) Except the case where $G = GL_n$, the above map in (ii) is not surjective.
In fact, the difference $\d = |\SN_G| - |W\wg| = 1$ if $G$ is of type $G_2, F_4$ or $E_8$, 
while $\d = 0$ if $G$ is of type $E_6$ or $E_7$. 
In the former case, the missing one is the pair $(x, \tau)$, where $x$ is the unique 
class such that $A_G(x) \simeq S_3, S_4, S_5$ (see Remark 3.1) and 
$\tau$ is the sign representation of $S_3, S_4, S_5$, respectively.  
\par
(ii) In the case of classical groups, such as $G = Sp_N$ or $SO_N$, 
$\d$ tends to $\infty$ if $N \to \infty$. 

\par\bigskip\noindent
{\bf \S 4.  The interpretation via perverse sheaves }
\addtocounter{section}{1}
\addtocounter{thm}{-3}
\par\medskip
In [L2] Lusztig reconstructed Springer representations of $W$ on $H^i(\SB_x, \Ql)$ 
in terms of the intersection cohomology. Based on Lusztig's 
construction, Borho-MacPherson [BM] reformulated Springer's theorem  (Theorem 3.2)
in the framework of the theory of perverse sheaves, which I will explain below.
\par
Let $\SN_G$ be the set of pairs $(\SO, \SE)$, where $\SO$ is a unipotent class 
in $G\uni$ and $\SE$ is a $G$-equivariant simple local system on $\SO$. If we choose 
$x \in \SO$, the stalk $\SE_x$ of $\SE$ at $x$ has a structure of a simple $A_G(x)$-module,
say $\tau \in A_G(x)\wg$, and $\tau$ characterizes a $G$-equivariant simple local system 
$\SE$, which we denote by $\SE_{\tau}$.  Then by the correspondence 
$(\SO, \SE) \lra (x, \tau)$, our set $\SN_G$ can be identified with the set 
$\SN_G$ appeared in Theorem 3.2 (ii). 
For each $(\SO, \SE) \in \SN_G$, we consider the intersection cohomology 
complex $\IC(\ol\SO, \SE)[\dim \SO]$, which gives a $G$-equivariant simple perverse 
sheaf on $G\uni$.  It is known by [L1] that the number of unipotent classes 
in $G\uni$ is finite. Thus by a general theory, we have
\par\medskip\noindent
(4.1) \ The set $\{ \IC(\ol\SO, \SE)[\dim \SO] \mid (\SO, \SE) \in \SN_G \}$
gives a complete set of isomorphism classes of $G$-equivariant simple local systems
on $G\uni$.  
\par\medskip
\begin{thm}[{Borho-MacPherson [BM]}]  
Let $\pi : \wt G\uni \to G\uni$ be as in Section 3, and consider 
the constant sheaf $\Ql$ on $\wt G\uni$.  
Then $\pi_*\Ql[\dim G\uni]$ is a semisimple perverse sheaf on $G\uni$, 
equipped with $W$-action, and is decomposed as 

\begin{equation*}
\pi_*\Ql[\dim G\uni] \simeq \bigoplus_{(\SO,\SE) \in \SN_G}
       V_{(\SO,\SE)}\otimes \IC(\ol\SO, \SE)[\dim \SO],
\end{equation*}
where $V_{(\SO,\SE)}$ is an irreducible $W$-module if it is non-zero. 
\end{thm}

Note that if $K = \pi_*\Ql[\dim G\uni]$ is a $G$-equivariant semisimple perverse
sheaf on $G\uni$, it is a direct sum of various $\IC(\ol\SO,\SE)[\dim \SO]$ 
by (4.1).  If $K$ is equipped with $W$-action, the multiplicity 
space $V_{(\SO,\SE)}$ has a structure of $W$-module.  The theorem asserts that
this $W$-module is irreducible.
Also note that the stalk $\SH^i_xK$ of the cohomology sheaf $\SH^iK$ at 
$x \in G\uni$ is isomorphic to $H^{i - \dim G\uni}(\SB_x, \Ql)$. Thus 
Theorem 3.2 is obtained as a corollary of Theorem 4.1.

\par\bigskip\noindent
{\bf \S 5. The generalized Springer correspondence }
\addtocounter{section}{1}
\addtocounter{thm}{-1}
\par\medskip
The map $W\wg \hra \SN_G$ is not necessarily surjective. It is an interesting 
problem to understand the pair $(\SO, \SE)$ which is not contained in the image 
of $W\wg$ (for example, the unique missing pair $(\SO, \SE)$ in the case of 
type $G_2, F_4$ or $E_8$ in Remarks 3.3).  
As an analogue of the Harish-Chandra theory of the representations of reductive 
groups, Lusztig extended this map to a bijection to $\SN_G$. His ingredients 
are as follows;
\par\medskip
\begin{itemize}
\item \ The notion of a cuspidal pair for $(\SO,\SE) \in \SN_G$, 
\item \ The notion of an induction $\ind_P^G$ for a parabolic subgroup $P$ of $G$ 
and its Levi subgroup $L$, 
\begin{equation*}
\ind_P^G : \{ \text{ $L$-equiv. perverse sheaves on $L$ } \}  
               \to \{ \text{ semisimple complexes on $G$ } \}. 
\end{equation*} 
\end{itemize}
\par\medskip\noindent
Note that $K$ is called a semisimple complex if it is a direct sum  of various 
$A[i]$, where $A[i]$ is a degree shift of a semisimple perverse sheaf $A$.
The induction functor $\ind_P^G$ is defined as follows.  We consider the diagram

\begin{equation*}
\begin{CD}
L  @<\a<<  \wh X_P  @>\psi>>  \wt X_P  @>\pi >>  G,
\end{CD}
\end{equation*}   
where
\begin{align*}
\wh X_P &= \{ (x, g) \in G \times G \mid g\iv xg \in P \}, \\
\wt X_P &= \{ (x, gP) \in G \times G/P \mid g\iv xg \in P \}, 
\end{align*}
and the maps are defined by 
\begin{equation*}
\a : (x, g) \mapsto \eta_P(g\iv xg), \quad 
\psi : (x, g) \mapsto (x, gP), \quad 
\pi : (x, gP) \mapsto x.
\end{equation*}
($\eta : P \to L \simeq P/U_P$ is the natural projection.) 
Here $\a$ is a smooth morphism with connected fibre, and $\psi$ 
is a principal $P$-bundle. 
Let $K$ be an $L$-equivariant perverse sheaf on $L$.  Then 
$\a^*K[a]$ is a $P$-equivariant perverse sheaf on $\wh X$ 
($a$ is the dimension of the fibre of $\a$), 
and there exists a unique perverse sheaf $\wt K$ on $\wt X$ such that
\begin{equation*}
\a^*K[a] \simeq \psi^*\wt K[b]
\end{equation*}
with $b = \dim P$. 
Since $\wt X_P$ is smooth, irreducible, and $\pi$ is proper,
$\pi_*\wt K$ is a semisimple complex on $G$ by the decomposition theorem 
of Deligne-Gabber.  We define  $\ind_P^GK = \pi_*\wt K$.   
Note that unless $P = B$, the map $\wt X_P  \to L$ is not defined directly,
so we need to consider $\wh X_P$ for defining $\wt K$.  
\par
Let $\ScS_G$ be the set of triples $(L, \SO_1, \SE_1)$, up to the natural action of $G$, 
where $L$ is a Levi subgroup of a parabolic subgroup of $G$, 
and $(\SO, \SE) \in \SN_L$ is a cuspidal pair. 
Put $\SW_L = N_G(L)/L$. In general, $\SW_L$ is not a Coxeter group, but 
in a very special situation that $\SN_L$ has a cuspidal pair, it turns out 
that $\SW_L$ is a Coxeter group.   
For $\xi = (L, \SO_1, \SE_1) \in \ScS_G$, put $K_{\xi} = \IC(\ol\SO_1, \SE_1)[\dim \SO_1]$.
Then $K_{\xi}$ is an $L$-equivariant perverse sheaf on $L$, and one can consider 
the complex $\ind_P^G K_{\xi}$ on $G$. 
The following result was proved by Lusztig.

\begin{thm}[{Lusztig [L3]}]  
\begin{enumerate}
\item
For $\xi = (L, \SO_1, \SE_1) \in \ScS_G$, $K_{\xi}$ is a semisimple perverse 
sheaf on $G$, equipped with $\SW_L$-action, and is decomposed as 

\begin{equation*}
\ind_P^GK_{\xi} \simeq \bigoplus_{(\SO, \SE) \in \SN_G }
                V^{\xi}_{(\SO, \SE)}\otimes \IC(\ol\SO, \SE)[\dim \SO],
\end{equation*}
where $V^{\xi}_{(\SO,\SE)}$ is an irreducible $\SW_L$-module if it is non-zero.
\item 
For any $\xi \in \ScS_G$, we have a bijection

\begin{equation*}
\SN_G^{(\xi)} := \{ (\SO, \SE) \in \SN_G \mid V^{\xi}_{(\SO,\SE)}\ne 0 \}
                     \isom \SW_L\wg. 
\end{equation*}
\item
We have a partition $\SN_G = \coprod_{\xi \in \ScS_G}\SN_G^{(\xi)}$.
In particular, there exists a natural bijection 

\begin{equation*}
\SN_G \isom \coprod_{\xi \in \ScS_G} \SW_L\wg,
\end{equation*}
which is called the {\bf generalized Springer correspondence. } 
\end{enumerate}
\end{thm}   
  
\remarks{5.2.} (i) \
In general, cuspidal pairs occur very rarely. 
Let $G = Sp_N$ or $SO_N$, and $\d_N$ the number of cuspidal pairs in $G$.
Then we have
\par\medskip\noindent
$\bullet$ \ $G = Sp_N$ : $\d_N = \begin{cases}
                           1 &\quad\text{ if $N = d(d-1)$ for some $d \in \BZ$ },  \\ 
                           0 &\quad\text{ otherwise. }                    
                                  \end{cases}$
\par\medskip\noindent
$\bullet$ \ $G = SO_N$ : $\d_N = \begin{cases}
                           1 &\quad\text{ if $N = d^2$ for some $d \in \BZ$ },  \\
                           0 &\quad\text{ otherwise. }
                                 \end{cases}$
\par\medskip
For example, the occurrence of the cuspidal pairs are given as

\begin{equation*}
\begin{aligned}
N &= 2, \ 6, \ 12, \ 20, \ 30, \dots &\quad &(G: Sp_N), \\
N &= 1, \ 4, \ 9, \ 16, \ 25, \dots &\quad  &(G : SO_N). 
\end{aligned}
\end{equation*}
\par 
(ii) \  If $\xi_0 = (L, \SO_1, \SE_1)$, where $L = T$ is a maximal torus, 
$\SO_1 = \{ 1 \}$ is the identity class, and $\SE_1$ is the constant 
sheaf $\Ql$, then $\ind_P^GK_{\xi_0}$ gives the original Springer correspondence
in Theorem 3.2 and Theorem 4.1. 
In this case, 
\par\medskip
$\bullet$ $1_W \lra \IC(\ol\SO_1, \Ql)$ : $\SO_1$ : the regular unipotent class in $G\uni$.  
\par\medskip
$\bullet$ $\ve_W \lra \IC(\ol\SO_0, \Ql)$ : $\SO_0 = \{ 1 \}$ ; the identity class.         

\par\bigskip\noindent
{\bf \S 6. Geometric realization of Kostka polynomials }
\addtocounter{section}{1}
\addtocounter{thm}{-2}
\par\medskip
In this section, we will discuss the relationship between Kostka polynomials and
the theory of the Springer correspondence in the case of $GL_n$.  In [L2], 
Lusztig gave a geometric realization of Kostka polynomials in terms of the 
intersection cohomology associated to unipotent classes in $GL_n$.  
First we review the definition of Kostka polynomials.
\par
Consider $x = (x_1, x_2, \dots,)$ infinitely many variables, and 
$t$ another parameter. Let $\vL(x) = \vL = \bigoplus_{n \ge 0}\vL^n$ be the ring 
of symmetric functions, where $\vL^n$ is the $n$-th homogeneous part, 
and consider $\vL[t] = \vL \otimes_{\BZ}\BZ[t]$, 
the ring over $\BZ[t]$.
For $\la \in \SP_n$, one can consider the Schur function $s_{\la}(x) \in \vL^n$
and the Hall-Littlewood function $P_{\la}(x;t) \in \vL^n[t] = \vL^n\otimes_{\BZ}\BZ[t]$.  
$\{ s_{\la}(x) \mid \la \in \SP_n \}$ and $\{ P_{\la}(x;t) \mid \la \in \SP_n \}$ 
give $\BZ[t]$-bases of $\vL^n[t]$. Thus one can define the {\bf Kostka polynomial} 
$K_{\la,\mu}(t) \in \BZ[t]$ by the condition that

\begin{equation*}
\tag{6.1}
s_{\la}(x) = \sum_{\mu \in \SP_n}K_{\la,\mu}(t)P_{\mu}(x;t).
\end{equation*}

Here we introduce two combinatorial objects.  
For $\la, \mu \in \SP_n$, write them as 
$\la = (\la_1,\dots, \la_k), \mu = (\mu_1, \cdots, \mu_k)$ with $\la_i, \mu_i \ge 0$.
We define a partial order 
$\la \le \mu$ on $\SP_n$, called the {\bf dominance order,} by the condition 
that, for any $1 \le j \le k$,
\begin{equation*}
\tag{6.2}
\sum_{i = 1}^j \la_i \le \sum_{i=1}^j\mu_i. 
\end{equation*}
On the other hand, for $\la \in \SP_n$, we define an {\bf $n$-function} 
$n : \SP_n \to \BZ_{\ge 0}$ by  
\begin{equation*}
\tag{6.3}
n(\la) = \sum_{i = 1}^k(i-1)\la_i. 
\end{equation*}
Then it is known that 
$K_{\la,\mu}(t)$ is a monic polynomial with 
$\deg K_{\la,\mu}(t) = n(\mu) - n(\la)$.  It follows that  
one can define a modified Kostka polynomial $\wt K_{\la,\mu}(t) \in \BZ[t]$ 
by 
\begin{equation*}
\tag{6.4}
\wt K_{\la,\mu}(t) = t^{n(\mu)}K_{\la,\mu}(t\iv).
\end{equation*}
\par
Let $\SO_{\la}$ be the unipotent class in $G\uni$ 
corresponding to $\la \in \SP_n$ for $G = GL_n$ as in Section 2. 
The following result shows that the closure relations of unipotent classes  
in $GL_n$ can be described in terms of the dominance order on $\SP_n$.
\begin{equation*}
\tag{6.5.}
\ol\SO_{\la} = \coprod_{\mu \le \la}\SO_{\mu}
\end{equation*}

The following result gives a geometric realization of 
Kostka polynomials. 
\begin{thm}[{Lusztig [L2]}]  
For $\la \in \SP_n$, put $A_{\la} = \IC(\ol\SO_{\la}, \Ql)$.  Then 
\begin{enumerate}
\item
$\SH^iA_{\la} = 0$ for odd $i$. 
\item
For $x \in \SO_{\mu} \subset \ol\SO_{\la}$, we have
\begin{equation*}
\tag{6.6}
\wt K_{\la,\mu}(t) = t^{n(\la)}\sum_{i \ge 0}(\dim \SH^{2i}_x A_{\la})t^i.
\end{equation*}
In particular, $K_{\la,\mu}(t) \in \BZ_{\ge 0}[t]$. 
\end{enumerate}
\end{thm}

In the case of $GL_n$, Borho-MacPherson's theorem gives a decomposition 
\begin{equation*}
\tag{6.7}
\pi_*\Ql[\dim G\uni] \simeq \bigoplus_{\la \in \SP_n}
         V_{\la} \otimes A_{\la}[\dim \SO_{\la}]. 
\end{equation*}
\par\noindent
By taking the stalk of the $i$-th cohomology on both sides of (6.7), 
and by comparing them   
with Theorem 6.1 (ii), we obtain the interpretation of Ksotka polynomials 
in terms of the Springer representations of $S_n$. 

\begin{cor}  
For $x \in \SO_{\mu}$, we have

\begin{equation*}
\wt K_{\la,\mu}(t) = t^{n(\la)}\sum_{i \ge 0}\lp H^{2i}(\SB_x, \Ql), V_{\la}\rp_{S_n}t^i,
\end{equation*}
where $\lp \ ,\ \rp_{\, S_n}$ is the inner product of characters of $S_n$. 
\end{cor}

\remark{6.3.}
In the case of the Springer correspondence, only Springer representations 
of the top cohomology $H^{2d_x}(\SB_x, \Ql)$ are involved.  While, 
for the description of the Kostka polynomials, all the Springer representations 
$H^{2i}(\SB_x, \Ql)$ are used. 

\par\bigskip\noindent
{\bf \S 7. The enhanced variety $GL(V) \times V$ }
\addtocounter{section}{1}
\addtocounter{thm}{-3}
\par\medskip
Later we discuss the theory of the Springer correspondence for 
the exotic symmetric space associated to symplectic groups. The enhanced 
variety is an extension of $GL_n$, and is also regarded as 
a degenerate form  of the exotic symmetric space, and 
was studied extensively by Achar-Henderson [AH] and Finkelberg-Ginzburg-Travkin [FGT]. 
So we will start from  the exposition on the enhanced variety. 
\par
Let $G = GL(V) \simeq GL_n$ with $\dim V = n$. 
We consider the direct product $X = G \times V$, on which $G$ acts as
$g : (x,v) \mapsto (gxg\iv, gv)$, where $gv$ is the natural action of $G$ on $V$. 
The role of $G\uni$ is played by $X\uni = G\uni \times V$, which is a $G$-invariant 
closed subset of $X$. The varieties $X_, X\uni$ are called the {\bf enhanced varieties}.
\par
Let $\SP_{n,2}$ be the set of double partitions $\Bla = (\la',\la'')$ such that 
$|\Bla| = |\la'| + |\la''| = n$. 
The following result is known by Achar-Henderson and Travkin. 
The finiteness of the $G$-orbits is crucial for later discussion. 

\begin{lem}  
Let $X\uni/\!\sim_G$ be the set of $G$-orbits in $X\uni$.  Then we have
\begin{equation*}
X\uni/\!\sim_G \, \simeq \SP_{n,2}.
\end{equation*} 
\end{lem}

The correspondence is given as follows; 
take $z = (x,v) \in G\uni \times V$.  Put $E^x = \{ y \in \End(v) \mid xy = yx \}$. 
Then $E^x$ is a subalgebra of $\End(V)$, and $x \in E^x$. 
Put $V_x = E^xv \subset V$.  Then $V_x$ is an $x$-stable subspace of $V$. 
Let $\la'$ be the Jordan type of $x|_{V_x}$ and $\la''$ the Jordan type of 
$x|_{V/V_x}$.   
Then $\Bla = (\la',\la'') \in \SP_{n,2}$, and the correspondence 
$z \mapsto \Bla$ gives a bijection $X\uni/\!\sim_G \, \isom \SP_{n,2}$. 
\par
We shall define an analogue of the dominance order on $\SP_{n,2}$ as follows;
for $\Bla = (\la',\la'') \in \SP_{n,2}$, write it as 
$\la' = (\la'_1, \dots, \la'_m), \la'' = (\la''_1, \dots, \la''_m)$ with common $m$, 
and define $c(\Bla) \in \BZ_{\ge 0}^{2m}$ by 

\begin{equation*}
\tag{7.1}
c(\Bla) = (\la_1', \la_1'', \la_2', \la_2'', \dots, \la'_m, \la''_m).
\end{equation*} 
Now define a dominance order $\Bla \le \Bmu$ on $\SP_{n,2}$ by the condition that
$c(\Bla) \le c(\Bmu)$ in $\BZ_{\ge 0}^{2m}$.  Note that the definition of the dominance
order on $\SP_n$ makes sense even for $\BZ_{\ge 0}^{2m}$. 
The following is known by [AH]. 

\begin{lem} 
For each $\Bla \in \SP_{n,2}$, we have $\ol\SO_{\Bla} = \coprod_{\Bmu \le \Bla}\SO_{\Bmu}$. 
\end{lem}

\remark{7.3.}
For $z = (x,v) \in G\uni \times V$, the stabilizer $Z_G(z)$ in $G$ is connected.
Hence the $G$-equivariant simple local system on a $G$-orbit $\SO$ in $X\uni$
is a constant sheaf $\Ql$, and we are in a quite similar situation as 
in the case of $GL_n$. 

\par\bigskip\noindent
{\bf \S 8. Springer correspondence for the enhanced variety }
\addtocounter{section}{1}
\addtocounter{thm}{-3}
\par\medskip
We fix a $B$-stable flag $M_1 \subset M_2 \subset \cdots \subset M_n = V$, 
with $\dim M_i = i$. For an integer $0 \le m \le n$, we define varieties

\begin{align*}
\tag{8.1}
\wt X_m &= \{ (x,v, gB) \in G\uni \times V \times \SB 
              \mid g\iv xg \in B, g\iv v \in M_m \}, \\
   X_m &= \bigcup_{g \in G}g(U \times M_m),  
\end{align*}
where $U$ is the unipotent radical of $B$, 
and define a map $\pi_{\Bm} : \wt X_m \to X_m$ by $(x,v, gB) \mapsto (x,v)$.
Then $\wt X_m$ is smooth, irreducible, and $\pi$ is proper. 
\par
For $0 \le m \le n$, put $\Bm = (m, n-m)$ and $S_{\Bm} = S_m \times S_{n-m}$.
Then $S_{\Bm}$ is a Weyl subgroup of $S_n$.  Put

\begin{equation*}
\tag{8.2}
\SP(\Bm) = \{ \Bla = (\la',\la'') \in \SP_{n,2} \mid |\la'| = m, |\la''| = n-m \}.
\end{equation*}
Then it is known  that

\begin{equation*}
\tag{8.3}
S_{\Bm}\wg \simeq \SP(\Bm), \quad V_{\Bla} = V_{\la'}\boxtimes V_{\la''}
                      \ \lra \ \Bla = (\la',\la'').
\end{equation*}

The following result is an analogue of Borho-MacPherson's theorem to the case of 
the enhanced variety. 
\begin{thm}[{[SS]}]   
Put $d_{\Bm} = \dim X_{\Bm}$.
\begin{enumerate}
\item
$(\pi_{\Bm})_*\Ql[d_{\Bm}]$ is a semisimple perverse sheaf on $X_{\Bm}$, 
equipped with $S_{\Bm}$-action, and is decomposed as 
\begin{equation*}
(\pi_{\Bm})_*\Ql[d_{\Bm}] \simeq \bigoplus_{\Bla \in \SP(\Bm)}
                      V_{\Bla} \otimes \IC(\ol\SO_{\Bla}, \Ql)[\dim \SO_{\Bla}].
\end{equation*}

\item
The Springer correspondence is given by 

\begin{equation*}
\coprod_{0 \le m \le n}(S_m \times S_{n-m})\wg \simeq 
         \coprod_{0 \le m \le n}\{ \SO_{\Bla} \mid \Bla \in \SP(\Bm) \}
              = X\uni/\!\sim_G.  
\end{equation*}
\end{enumerate}
\end{thm}

\par\bigskip\noindent
{\bf \S 9. Double Kostka polynomials }
\addtocounter{section}{1}
\addtocounter{thm}{-1}
\par\medskip
Kostka polynomials $K_{\la,\mu}(t)$ ($\la,\mu \in \SP_n$) was generalized 
to double Kostka polynomials $K_{\Bla,\Bmu}(t)$ ($\Bla,\Bmu \in \SP_{n,2}$). 
(See [S1, S2], where they are defined for any $r$-partitions. 
Here we only consider the special case where $r = 2$). 
Recall the discussion in Section 6. Here we prepare two types of 
variables $x' = (x_1',x_2', \dots)$, $x'' = (x''_1, x''_2, \dots)$
and consider $\Xi = \bigoplus_{n \ge 0}\vL(x')\otimes \vL(x'')$, the ring
of symmetric functions with respect to $x = (x',x'')$. 
Put $\Xi[t] = \Xi\otimes_{\BZ}\BZ[t]$. 
For $\Bla = (\la',\la'') \in \SP_{n,2}$, one can define the Schur function 
$s_{\Bla}(x) = s_{\la'}(x')s_{\la''}(x'')$ and Hall-Littlewood functions 
$P_{\Bla}(x;t)$.  (Note that the definition of $P_{\Bla}(x;t)$ is 
rather complicated.)  Then $\{ s_{\Bla}(x) \mid \Bla \in \SP_{n,2}\}$, 
$\{ P_{\Bla}(x;t) \mid \Bla \in \SP_{n,2} \}$ give two $\BZ[t]$-bases 
of $\Xi^n[t] = \Xi^n\otimes_{\BZ}\BZ[t]$. We define the double Kostka polynomial
$P_{\Bla,\Bmu}(t) \in \BZ[t]$ by 
\begin{equation*}
\tag{9.1}
s_{\Bla}(x) = \sum_{\Bmu \in \SP_{n,2}}K_{\Bla,\Bmu}(t)P_{\Bmu}(x;t). 
\end{equation*}
\par
For $\Bla = (\la', \la'') \in \SP_{n,2}$, define an $a$-function 
$a : \SP_{n,2} \to \BZ_{\ge 0}$ by

\begin{equation*}
\tag{9.2}
a(\Bla) = 2(n(\la') + n(\la'')) + |\la''|.
\end{equation*}
The $a$-function has a role of $n$-function in the case of $\SP_{n,2}$. 
It is known that $K_{\Bla,\Bmu}(t)$ is a monic of degree $a(\Bmu) - a(\Bla)$. 
As an analogue of (6.4), we define a modified double Kostka polynomial 
$\wt K_{\Bla,\Bmu}(t)$ by 

\begin{equation*}
\tag{9.3}
\wt K_{\Bla,\Bmu}(t) = t^{a(\Bmu)}K_{\Bla,\Bmu}(t\iv).
\end{equation*}

\par
The following result was proved by [AH], which is an analogue of 
Theorem 6.1, and gives a geometric 
realization of double Kostka polynomials. 

\begin{thm}[{[AH]}]   
For $\Bla \in \SP_{n,2}$, put $A_{\Bla} = \IC(\ol\SO_{\Bla}, \Ql)$. 
Then we have
\begin{enumerate}
\item
$\SH^iA_{\Bla} = 0$ for odd $i$.
\item
Take $\Bla,\Bmu \in \SP_{n,2}$ such that $\Bmu \le \Bla$. 
For $z \in \SO_{\Bmu} \subset \ol\SO_{\Bla}$, we have

\begin{equation*}
\tag{9.4}
\wt K_{\Bla,\Bmu}(t) = t^{a(\Bla)}\sum_{i \ge 0}(\dim \SH^{2i}_zA_{\Bla})t^{2i}.
\end{equation*} 
\end{enumerate}
\end{thm}

\remark{9.2}  
If we compare the formula (9.4) with the formula (6.6) in Theorem 6.1, 
we notice the discrepancy of the relations between $t^i$ and cohomology sheaves.
In fact, in (6.6) $t^i$ corresponds to $\SH^{2i}_xA_{\la}$, while in (9.4)
$t^{2i}$ corresponds to $\SH^{2i}_zA_{\Bla}$. Later, we also consider the geometric
realization of double Kostka polynomials by means of exotic symmetric space.
In that case, this discrepancy is removed (see Theorem 13.4). 

\par\bigskip\noindent
{\bf \S 10. Symmetric spaces in algebraic setting}
\addtocounter{section}{1}
\addtocounter{thm}{-2}
\par\medskip
Historically, the symmetric space $G/K$, where $G$ is a connected Lie group and 
$K$ is a compact subgroup such that $(G^{\th})^0 \subset K \subset G^{\th}$
for an involutive automorphism $\th$ of $G$, has been studied extensively 
from a point of view of the Riemannian geometry.  However, here we are interested 
in its algebraic structure, namely, we consider $G, K$ as algebraic groups 
and $G/K$ as an algebraic variety over $\Bk$ of any characteristic.  
The basis of the algebraic study of symmetric spaces was achieved by Vust [V] 
and Richardson [R], which will be summarized below. 
\par
Let $G$ be a connected reductive group over $\Bk$ with $\ch \Bk \ge 0$, 
and $\th : G \to G$ 
an automorphism such that $\th^2 = 1$. Consider the fixed point subgroup 
$G^{\th} = \{ g \in G \mid \th(g) = g \}$ of $G$, and put $H = (G^{\th})^0$. 
(Here we use the notation $H$ rather than $K$, since we often use $K$ as 
complexes of sheaves.) Put
\begin{equation*}
G^{\io\th} = \{ g \in G \mid \th(g) = g\iv \},
\end{equation*}   
where $\io : G \to G, g \mapsto g\iv$ is an anti-automorphism, 
so we regard $G^{\io\th}$ as the subset of $G$ consisting of $\io\th$-fixed 
elements. Put $G^{\io\th}_0 = \{ g\th(g)\iv \mid g \in G \}$. 
It is known that $G^{\io\th}_0$ is a connected component of $G^{\io\th}$, 
and there exists an isomorphism 

\begin{equation*}
G^{\io\th}_0 \isom G/G^{\th}, \quad g\th(g)\iv \lra gG^{\th}.
\end{equation*}

$G^{\th}$ acts on $G^{\io\th}$ and $G^{\io\th}_0$ by the conjugation action, 
which corresponds to the left multiplication of $G^{\th}$ on $G/G^{\th}$ via 
the above isomorphism. 
Thus one can identify $G^{\io\th}_0$ with the symmetric space $G/G^{\th}$, and 
from now on, we consider $G_0^{\io\th}$ as a symmetric space in an algebraic setting. 

\remark{10.1} 
Assume that $p \ne 2$.  Let $\th : G \to G$ be as above.  
Let $\Fg$ be the Lie algebra of $G$. 
Then $\th$ induces an involutive automorphism on $\Fg$, which we denote by 
the same symbol as $\th : \Fg \to \Fg$.  We have the decomposition 
\begin{equation*}
\Fg \simeq \Fg^{\th} \oplus \Fg^{-\th},
\end{equation*}
where $\Fg^{\pm \th}$ are $\pm 1$ eigenspaces of $\th$. We note that 
$\Fg^{\th} \simeq \Lie G^{\th}$ and $\Fg^{-\th}$ is a $G^{\th}$-stable subspace
of $\Fg$ such that $\Fg^{-\th} \simeq T_e(G^{\io\th})$. 
$\Fg^{-\th}$ is usually referred as the symmetric space with respect to $G^{\th}$.
Thus $G^{\io\th}$ is regarded as a global analogue of $\Fg^{-\th}$. 
\par\medskip
Let $G^{\io\th}\uni = G^{\io\th}_0 \cap G\uni$ be the set of unipotent elements in 
$G^{\io\th}_0$.  Then $G^{\io\th}\uni$ is an $H$-stable closed subset of $G^{\io\th}$, 
which has a role of the unipotent variety $G\uni$ in the case of symmetric 
spaces.
The following important result was proved by Richardson [R].
In fact, Richadson reduced the problem to a similar problem for reductive groups, 
which certainly holds by Lusztig [L1]. 

\begin{prop}  
Assume that $p \ne 2$.  Then the number of $H$-orbits in $G^{\io\th}_0$ is finite. 
\end{prop}   

\remark{10.3.}
In the Lie algebra case, put $\Fg^{-\th}\nil = \Fg^{-\th} \cap \Fg\nil$, 
where $\Fg\nil$ is the set of nilpotent elements in $\Fg$.
Then $\Fg^{-\th}\nil$ is an $H$-stable closed subset of $\Fg^{-\th}$. 
The finiteness property also holds in the Lie algebra case.  
But this does not hold for $G^{\io\th}\uni$ if $p = 2$. 

\par\bigskip\noindent
{\bf \S 11. Symmetric spaces associated to classical groups }
\par\medskip
From now on, we concentrate ourselves to the special type of symmetric 
spaces, namely the symmetric spaces associated to classical groups.   
\par
Let $V$ be an $N$-dimensional vector space over $\Bk$ with $\ch\Bk \ne 2$, 
and put $G = GL_N \simeq GL(V)$. 
We define an involutive automorphism $\th : G \to G$ by 
$\th(g) = J\iv ({}^tg\iv)J$, where $J = J_1$ or $J_2$, 

\begin{align*}
J_1 &= \begin{pmatrix}
         1  &  0    &  0   \\
         0  &  0    &  1_n  \\
         0  &  1_n    &  0
       \end{pmatrix}
          \quad \text{ if $N = 2n+1$, }  \quad
J_1 =  \begin{pmatrix}
            0      &   1_n  \\
            1_n    &   0
       \end{pmatrix}
          \quad\text{ if $N = 2n$, }   \\
J_2 &=  \begin{pmatrix}
             0    &   1_n  \\
            -1_n  &   0    
        \end{pmatrix}    
           \quad\text{ for $N = 2n$. }
\end{align*}
If $J = J_1$, then $G^{\th} = O_N$ and $H = SO_N$.  While if 
$J = J_2$, then $G^{\th} = H = Sp_N$. 
In both cases, the following identity holds. 

\begin{equation*}
\tag{11.1}
G^{\io\th}_0 = G^{\io\th}.  
\end{equation*}

For later applications, we consider a generalization of the symmetric 
space $G^{\io\th}$ of the following type; 
for an integer $r \ge 1$, consider the direct product 
$G^{\io\th}\uni \times V^{r-1}$, on which $H$ acts diagonally. 
$G^{\io\th}\uni \times V^{r-1}$ is called the {\bf exotic symmetric space} 
of level $r$.  In the following discussion, we are interested in extending 
the theory of the Springer correspondence to the case of exotic symmetric spaces.
\par
Let $T \subset B$ be the pair of $\th$-stable maximal torus and $\th$-stable 
Borel subgroup of $G$. Then the unipotent radical $U$ of $B$ is $\th$-stable. 
Put $B_H = (B^{\th})^0$ and $T_H = (T^{\th})^0$. 
Then $B_H$ is a Borel subgroup of $H$ and $T_H$ is a maximal torus of $H$ contained 
in $B_H$.  We consider the flag variety $\SB^H = H/B_H$ of $H$. 
Let $M_1 \subset \cdots \subset M_n$ be the isotropic flag in $V$ such that  
the stabilizer of $(M_i)_i$ in $H$ is $B_H$.
Consider the varieties

\begin{align*}
\wt\SX &= \{ (x, \Bv, gB_H) \in G^{\io\th}\uni \times V^{r-1} \times \SB^H
              \mid g\iv xg \in B^{\io\th}, g\iv\Bv \in M_n^{r-1} \}, \\
   \SX &= \bigcup_{g \in H}g(U^{\io\th} \times M_n^{r-1}).
\end{align*}  
and define a map $\pi : \wt\SX \to \SX$ by 
$(x,\Bv, gB_H) \mapsto (x,\Bv)$. 
\par
The map $\pi : \wt\SX \to \SX$ is an analogue of the Springer resolution 
$\pi : \wt G\uni \to G\uni$, and we want to consider an analogue of the 
Borho-MacPherson's theorem. But in order to apply the previous discussion, 
we need to verify two crucial properties, namely (i) the map $\pi$ gives 
a resolution of singularities, (ii)  the number of $H$-orbits is finite. 
(Actually, (i) is too strong. For the theory of the Springer correspondence, 
enough to show that $\pi$ is ``semi-small'', which implies that 
(i$'$) $\dim \wt\SX = \dim \SX$).
Thus we will verify those two properties (i$'$) and (ii).   Put 
$\d = \dim \wt\SX - \dim \SX$.  The following holds. 

\par\bigskip
$\bullet$ \ $H = Sp_N$ 
\par
\vspace{-1cm}
\begin{table}[h]
\hspace*{-2cm}
\begin{center}
\begin{tabular}{c|c|c|c}
   $r$ \ &  $r = 1$  &  $r = 2$  &  $r \ge 3$   \\  
\hline 
 $\d$ \    &   $\d > 0$   &  $\d = 0$  &  $\d = 0$          \\
\hline  
 number of orbits &  $< \infty$ &  $< \infty$   &   $\infty$ 
\end{tabular}
\end{center}
\end{table}

$\bullet$ \ $H = SO_N$

\vspace{-1cm}
\begin{table}[h]
\hspace*{-2cm}
\begin{center}
\begin{tabular}{c|c|c}
   $r$ \ &  $r = 1$  &  $r \ge 2$  \\  
\hline 
 $\d$ \    &   $\d = 0$   &  $\d = 0$  \\
\hline  
 number of orbits &  $< \infty$ &  $\infty$ 
\end{tabular}
\end{center}
\end{table}

\par
From those tables, we see that the most suitable situation for the 
Springer correspondence is that $\SX = G^{\io\th}\uni \times V$ for 
the case $H = Sp_N$, and $\SX = G^{\io\th}\uni$ for the case $H = SO_N$.  
This gives a reason why considering the exotic symmetric space, rather 
then the symmetric space itself, is important.  In the case where $H = Sp_N$, 
it is more natural to consider $G^{\io\th}\uni \times V$ than $G^{\io\th}\uni$.
In the following sections, we consider those standard cases, separately. 
However, we consider the other cases also, since some modified theory of 
the Springer correspondence holds in those cases, and they have own interests.

\par\bigskip\noindent
{\bf \S 12.  Unipotent orbits in $H$ and in $G^{\io\th}\uni$ }
\addtocounter{section}{2}
\addtocounter{thm}{-3}
\par\medskip
For later applications, in this section we describe the $H$-orbits 
in $G^{\io\th}\uni$ in connection with the unipotent classes in $G^{\th}\uni$. 
\par
In the following discussion, we only consider $H = SO_N$ ($N$ : odd) or $H = Sp_N$, 
namely $H$ is of type $B_n$ or $C_n$. 
Similar results hold in the case where $H$ is of type $D_n$, but since 
the description becomes more complicated, we omit this, for simplicity. 
\par
Let $\la  = (1^{m_1}, 2^{m_2}, \dots ) \in \SP_N$ be a partition of $N$. 
Let $\OO_{\la}$ be the unipotent class in $G = GL_N$ corresponding to $\la$.
The following result is well-known.
\begin{prop}[{\bf group case}]  
Under the notation above, we have 
\begin{enumerate}
\item
$\OO_{\la} \cap H \ne \emptyset$ $\Longleftrightarrow$ 
$m_i$ is even for odd $i$ $($ resp. for even $i$ $)$ \\
if $H = Sp_N$ $($ resp. $H = SO_N$ $)$.

\item
If $\OO_{\la} \cap H \ne \emptyset$, then 
$\SO_{\la} = \OO_{\la} \cap H$ is a single class in $H$.

\item
For $x \in \SO_{\la}$, 
\par\vspace{-8mm}
\begin{equation*}
Z_{G^{\th}}(x) \simeq \begin{cases}
           \bigl(\prod_{i : \odd }Sp_{m_i}
             \times  \prod_{i : \even }O_{m_i}\bigr) \ltimes U_1
                    &\quad\text{ $($ $Sp_N$-case $)$, } \\
            \bigl(\prod_{i : \even }Sp_{m_i}
               \times  \prod_{i : \odd }O_{m_i}\bigr) \ltimes U_1
                    &\quad\text{ $($ $SO_N$-case $)$, } 
                          \end{cases}
\end{equation*}
where $U_1$ : connected unipotent normal group of $Z_{G^{\th}}(x)$. 
\end{enumerate}
\end{prop}

\par
The case of $H$-orbits in $G^{\io\th}\uni$ is given as follows.

\begin{prop}[{\bf symmetric space case}]
Let $H = Sp_N$.  
\begin{enumerate}
\item 
$\OO_{\la} \cap G^{\io\th} \ne \emptyset$ $\Longleftrightarrow$ 
$m_i$ is even for all $i$.
\item 
If $\OO_{\la} \cap G^{\io\th} \ne \emptyset$, then 
$\SO_{\la} = \OO_{\la} \cap G^{\io\th}$ is a single $H$-orbit.
\item 
For $x \in \SO_{\la}$, 
\begin{equation*}
Z_{G^{\th}}(x) \simeq  \prod_iSp_{m_i} \ltimes U_1.
\end{equation*}
\end{enumerate}
\end{prop}

\begin{prop}[{\bf symmetric space case}]
Let $H = SO_N$ $($ $N$ : odd $)$.
\begin{enumerate}
\item
$\OO_{\la} \cap G^{\io\th}$ is always non-empty.
\item
$\SO_{\la} = \OO_{\la} \cap G^{\io\th}$ is a single $H$-orbit. 
\item
For $x \in \SO_{\la}$, 
\begin{equation*}
Z_{G^{\th}}(x) \simeq \prod_iO_{m_i} \ltimes U_1,
\end{equation*}
\end{enumerate} 
\end{prop}

\remark{12.4}
By comparing Proposition 12.2 and Propositions 12.3, 12.4, we find 
very interesting phenomena.   
In the group case, $Z_{G^{\th}}(x)$ involves the subgroups of
type $Sp$ and  $SO$.  Hence in the study of unipotent classes in $H$, the case 
of $Sp$ and of $SO$ cannot be separated. 
While in the case of symmetric spaces, the structure of $Z_{G^{\th}}(x)$ is 
completely separated, namely in the case of $H = Sp_N$, $Z_{G^{\th}}(x)$ involves 
only subgroups of type $Sp$, and similarly for $H = SO_N$. 
All the properties of $Z_{G^{\th}}(x)$ in the case of symmetric spaces are 
inherited from  similar properties in the group case, but they produce two extreme
phenomena to the opposite directions.  

\par\medskip
As corollaries of above results, we can determine the structure of 
the component group $A_H(x)$  as follows.   

\begin{cor}[{\bf group case}]  
Assume $x \in \SO_{\la} \subset H\uni$.  
Put $A_H(x) = Z_H(x)/Z^0_H(x)$. 
Then
\begin{equation*}
A_H(x) \simeq (\BZ/2\BZ)^{a(\la)}, 
\end{equation*}

where
\begin{equation*}
a(\la) = \begin{cases}
           \sharp\{ i \mid i : \text{ even, } m_i \ne 0 \}
              &\quad\text{ if $H = Sp_N$, } \\
           \sharp\{ i \mid i : \text{ odd, } m_i \ne 0 \} - 1 
               \ (\text { or } 0 )
              &\quad\text{ if $H = SO_N$. }  
         \end{cases}
\end{equation*}
\end{cor}

\begin{cor}[{\bf symmetric space case}]
Assume $x \in \SO_{\la} \subset G^{\io\th}\uni$.  
Put $A_H(x) = Z_H(x)/Z_H^0(x)$.  
Then 
\begin{equation*}
A_H(x) \simeq (\BZ/2\BZ)^{b(\la)},
\end{equation*}

where 
\begin{equation*}
b(\la) = \begin{cases}
            0   &\quad\text{ if $H = Sp_N$, } \\
            \sharp\{ i \mid m_i \ne 0 \} - 1 \ (\text{ or } 0 )
                &\quad\text{ if $H = SO_N$. } 
         \end{cases}
\end{equation*}
\end{cor}

\par\bigskip\noindent
{\bf \S 13. Exotic symmetric space associated to symplectic groups }
\addtocounter{section}{1}
\addtocounter{thm}{-6}
\par\medskip
We now consider the case of exotic symmetric space associated to symplectic
groups.  So assume that $H = Sp_N$ and $\SX = G^{\io\th}\uni \times V$. 
The fact that the set of $H$-orbits in $\SX$ is parametrized by $\SP_{n,2}$ 
was first proved by Kato [K1].  The following is a reformulation of Kato's result 
due to [AH].
\par
Let $M_n$ be the maximal isotropic subspace in $V$, stable by $B_H$. 
There exists a $\th$-stable Levi subgroup $L \subset G$ such that 
$L^{\th} \simeq GL(M_n)$.  Then we have embeddings 

\begin{equation*}
\tag{13.1}
L^{\io\th}\uni \times M_n \subset G^{\io\th}\uni \times V  \subset G\uni \times V.
\end{equation*} 
Note that the left hand side and the right hand side of (13.1) are enhanced 
varieties discussed in Section 7; 
$GL(M_n)$ acts on $L^{\io\th}\uni \times M_n$ diagonally, and their
orbits are parametrized by $\SP_n$ by Lemma 7.1.  In turn, $G$ acts on 
$G\uni \times V$ diagonally, and their orbits are parametrized by $\SP_{2n}$ by 
Lemma 7.1.  
We have

\begin{lem}[{[K1], [AH]}]  
Let $\Bla \in \SP_{n,2}$, and $\Bla \cup \Bla \in \SP_{2n,2}$.
Let $\SO^L_{\Bla}$ be the $GL(M_n)$-orbit corresponding to  $\Bla$, and 
$\OO_{\Bla \cup \Bla}$ the $G$-orbit corresponding to $\Bla \cup \Bla$. 
Then $\SO_{\Bla}$ is the unique $H$-orbit in $\SX$ 
such that
\begin{equation*}
\SO_{\Bla}^L \subset \SO_{\Bla} \subset \OO_{\Bla\cup\Bla}.
\end{equation*}
\end{lem}

The analogue of Lemma 7.2 and Remark 7.3 also hold for $G^{\io\th}\uni \times V$.

\begin{prop}  
Let $\SO_{\Bla}$ be an $H$-orbit in $G^{\io\th}\uni \times V$ corresponding to
$\Bla \in \SP_{n,2}$.
\begin{enumerate}
\item
For $z = (x,v) \in G^{\io\th}\uni \times V$, $Z_H(z)$ is connected.
\item
$\ol\SO_{\Bla} = \coprod_{\Bmu \le \Bla}\SO_{\Bmu}$, namely the closure 
relations are described by the dominance order on $\SP_{n,2}$. 
\end{enumerate} 
\end{prop}

In this case, the map $\pi : \wt\SX \to \SX$ in  Section 11 is given as follows;

\begin{align*}
\wt\SX &= \{ (x, v, gB_H) \in G^{\io\th}\uni \times V \times \SB^H
             \mid g\iv xg \in B^{\io\th}, g\iv v \in M_n \}, \\
   \SX &= \bigcup_{g \in H}g(U^{\io\th} \times V) = G^{\io\th}\uni \times V,
\end{align*}
and $\pi : (x, v, gB_H) \mapsto (x,v)$. 
\par
Let $W_n$ be the Weyl group of $H = Sp_N$, namely the Weyl group of 
type $C_n$.  It is known that $W_n\wg \simeq \SP_{n,2}$, we denote by 
$\wt V_{\Bla}$ irreducible representation of $W_n$ corresponding to $\Bla \in \SP_{n,2}$. 
(Note the difference from  $V_{\Bla}$, which was defined in (8.3).) 
The following theorem was first proved by Kato [K1] by applying the Ginzburg 
theory of Hecke algebras, and then reproved by [SS, I] in the framework of 
Lusztig's theory of the generalized Springer correspondence. 

\begin{thm}  
$\pi_*\Ql[\dim \SX]$ is a semisimple perverse sheaf on $\SX$, equipped with $W_n$-action, 
and is decomposed as 

\begin{equation*}
\pi_*\Ql[\dim \SX] \simeq \bigoplus_{\Bla \in \SP_{n,2}}
                            \wt V_{\Bla}\otimes \IC(\ol\SO_{\Bla},\Ql)[\dim \SO_{\Bla}].
\end{equation*}
\end{thm} 

The geometric realization of double Kostka polynomials in terms of $G^{\io\th}\uni \times V$
was conjectured in [AH], and was proved by [K3], [SS, II], independently.
In [SS,II], it is proved by making use of the discussion on character sheaves due 
to Lusztig [L4].   

\begin{thm} 
For $\Bla \in \SP_{n,2}$, put $A_{\Bla} = \IC(\ol\SO_{\Bla}, \Ql)$.
Then we have
\begin{enumerate}
\item 
$\SH^iA_{\Bla} = 0$ unless $i \equiv 0 \pmod 4$.
\item
For $z \in \SO_{\Bmu} \subset \ol\SO_{\Bla}$, 
\begin{equation*}
\tag{13.2}
\wt K_{\Bla,\Bmu}(t) = t^{a(\Bla)}\sum_{i \ge 0}(\dim \SH^{4i}_zA_{\Bla})t^{2i}.
\end{equation*}
\end{enumerate}
\end{thm}

Compare the formula (13.2) with the formula (9.4) in the enhanced case.  
The correspondence $\SH^{4i} \lra t^{2i}$ in (13.2)
is more natural than the correspondence $\SH^{2i} \lra t^{2i}$ in (9.4).  
Note that the modulo 4 vanishing of the cohomology of $A_{\Bla}$ was 
first noticed by Grojnowski [Gr].
\par
For $z = (x,v) \in G^{\io\th}\uni \times V$, consider the Springer fibre 
\begin{equation*}
\pi\iv(z) \simeq \SB^H_z = \{ gB_H \in \SB^H \mid g\iv xg \in B^{\io\th}, g\iv v \in M_n \}. 
\end{equation*}
Then $H^i(\SB^H_z, \Ql)$ has a structure of $W_n$-module, which is an analogue 
of the Springer representations. By using the Springer representation of $W_n$, 
we obtain an expression of double Kostka polynomials, which is an analogue 
of Corollary 6.2. 

\begin{cor}  
Assume that $z \in \SO_{\Bmu}$.  Then 
$H^i(\SB^H_z, \Ql) = 0$ unless $i \equiv 0 \pmod 4$, and we have
\begin{equation*}
\wt K_{\Bla, \Bmu}(t) = t^{a(\Bla)}\sum_{i \ge 0}
          \lp H^{4i}(\SB^H_z, \Ql), \wt V_{\Bla}\rp_{\,W_n} t^{2i}.
\end{equation*}
\end{cor}

\par\bigskip\bigskip\noindent
{\bf \S 14.  Springer correspondence for $G^{\io\th}\uni$, the case $H = Sp_N$ }
\addtocounter{section}{1}
\addtocounter{thm}{-5}
\par\medskip
As pointed out in Section 11, $G^{\io\th}\uni$ for $H = Sp_N$ does not satisfy 
the condition explained there, so the Springer correspondence does not hold 
in the strict sense.  However some modified theory of the Springer correspondence
still holds for $G^{\io\th}\uni$, as discussed in Henderson [H], which I will explain below. 
\par
First we prepare some notation which is common for $H = Sp_N$ and $SO_N$. 
Assume that $H = Sp_N$ or $SO_N$ ($N$ : odd), and consider $G^{\io\th}\uni$.  
Let $\SN_{G^{\io\th}}$ be the set of pairs $(\SO, \SE)$, where $\SO$ is an $H$-orbit in
$G^{\io\th}\uni$, and $\SE$ is an $H$-equivariant simple local system on $\SO$. 
Thus as in Section 4, $\SN_{G^{\io\th}}$ can be expressed as follows;
\begin{equation*}
\tag{14.1}
\SN_{G^{\io\th}} \simeq \{ (x,\tau) \mid x \in G^{\io\th}\uni/\!\sim_H, \tau \in A_H(x)\wg \}. 
\end{equation*}
  
We consider the variety 

\begin{equation*}
\tag{14.2}
\wt G^{\io\th}\uni = \{ (x, gB_H) \in G^{\io\th}\uni \times \SB^H
                            \mid g\iv xg \in B^{\io\th} \}
\end{equation*}
and define a map $\pi : \wt G^{\io\th}\uni \to G^{\io\th}\uni$ by 
$(x, gB_H) \mapsto x$.
Then $\pi$ is proper, surjective, and $\wt G^{\io\th}\uni$ is smooth, 
irreducible. If we put $\d = \dim \wt G^{\io\th}\uni - \dim G^{\io\th}\uni$, 
$\d = 0$ for $H = SO_N$, while $\d > 0$ for $H = Sp_N$.  
In fact, in the latter case, $\pi$ gives rise to a $(\BP_1)^n$-bundle on its open dense
part, and $\d = n$. 
Let $W_n = N_H(T_H)/T_H$ be the Weyl group of $H$, thus $W_n$ is 
of type $C_n$ or $B_n$ according to the case where $H = Sp_N$ or $SO_N$. 
We consider $S_n$ as a subgroup of $W_n$. 
\par
Now assume that $H = Sp_N$.  
By Proposition 12.2, the set $G^{\io\th}\uni/\!\sim_H$ is parametrized by 
$\SP_n$,  under the correspondence 
\begin{equation*}
\tag{14.3}
\SO_{2\la} = \OO_{2\la} \cap G^{\io\th}\uni \lra \la \in \SP_n,
\end{equation*}   
and by Corollary 12.6, $Z_H(x)$ is always connected for $x \in G^{\io\th}\uni$. 
Hence the local system $\SE$ on $\SO$ is the constant sheaf $\Ql$, 
and we have  

\begin{equation*}
\tag{14.4}
\SN_{G^{\io\th}} =  \{ (\SO, \Ql) \} \simeq \SP_n. 
\end{equation*}  

The following result was first proved by [H] for $\Fg^{-\th}\nil$, 
without the explicit correspondence, and was 
proved for $G^{\io\th}\uni$ in [SS, I] in the following form.

\begin{thm}  
Let $H = Sp_N$. Then $\pi_*\Ql[\dim G^{\io\th}\uni]$ is a semisimple 
complex, equipped with $S_n$-action, and is decomposed as 

\begin{equation*}
\pi_*\Ql[\dim G^{\io\th}\uni] \simeq H^{\bullet}(\BP_1^n)\otimes 
    \bigoplus_{\la \in \SP_n}V_{\la} \otimes \IC(\ol\SO_{2\la}, \Ql)[\dim \SO_{2\la}], 
\end{equation*} 
where $H^{\bullet}(\BP_1^n) = \bigoplus_{i \ge 0}H^{2i}(\BP_1^n, \Ql)$ 
is a complex of vector spaces, and 
$V_{\la}$ is the irreducible $S_n$-module corresponding to $\la \in \SP_n$.  
\end{thm}
\par\noindent

Note that the occurrence of the factor $H^{\bullet}(\BP_1^n)$ 
depends on the fact that $\pi$ is generically a $\BP_1^n$-bundle. 

\remark{14.2.} 
If we consider the symmetric space $G^{\io\th}\uni$ of general type, 
one can expect that a similar phenomenon as in the case of $H = Sp_N$ occurs.  
The case where $H = SO_N$ seems to be rather special, and we could 
not find other examples such that the Springer correspondence holds 
in the strict sense.  Thus in order to discuss the (generalized) 
Springer correspondence for the general case, 
we nned to consider a correction factor such as 
$H^{\bullet}(\BP_1^n)$ in the case of $Sp_N$. 
In the Lie algebra case, Lusztig-Yun [LY] discusses a related problem.   

\par\bigskip\noindent
{\bf \S 15.  Springer correspondence for $G^{\io\th}\uni$, the case $H = SO_{2n+1}$ }
\addtocounter{section}{1}
\addtocounter{thm}{-1}
\par\medskip
We now assume that $H = SO_{2n+1}$.  Then by Proposition 12.3, 
the set of $H$-orbits in $G^{\io\th}\uni$ is parametrized by $\SP_{2n+1}$;
\begin{equation*}
\tag{15.1}
\SO_{\la} = \OO_{\la} \cap G^{\io\th}\uni \lra \la \in \SP_{2n+1}.
\end{equation*}

We define a map $\vG : \SP_n \to \SP_{2n+1}$ by 

\begin{equation*}
\tag{15.2}
\vG : \mu = (\mu_1, \dots, \mu_k) \mapsto (2\mu_1 +1, 2\mu_2, \dots, 2\mu_k).
\end{equation*}
The following result gives the Springer correspondence for $G^{\io\th}\uni$. 
It was first proved by Chen-Vilonen-Xue [CVX] 
for $\Fg^{-\th}\nil$ with $\Bk = \BC$. The group case $G^{\io\th}\uni$ 
is due to [SY].

\begin{thm}  
Assume that $H = SO_{2n+1}$. 
Then $\pi_*\Ql[\dim G^{\io\th}\uni]$ is a semisimple perverse sheaf 
on $G^{\io\th}\uni$, equipped with $S_n$-action, and is decomposed as

\begin{equation*}
\tag{15.3}
\pi_*\Ql[\dim G^{\io\th}\uni] \simeq \bigoplus_{\mu \in \SP_n}
           V_{\mu} \otimes \IC(\ol\SO_{\vG(\mu)}, \Ql)[\dim \SO_{\vG(\mu)}],
\end{equation*}
where $V_{\mu}$ is the irreducible $S_n$-module corresponding to $\mu \in \SP_n$. 
\end{thm}

\remarks{15.2.}
(i) \ The formula (15.3) looks very similar to the formula 
in the case of $GL_n$.  But the pattern of the Springer correspondence 
is quite different from  that of $GL_n$ (see Remarks 5.2 (ii)).   
In fact, the Springer correspondence for $G^{\io\th}\uni$ is given as 
follows;

\begin{itemize}
\item  $1_{S_n} = V_{\mu}$ with $\mu = (n) \lra \IC(\ol\SO_{\vG(\mu)}, \Ql)$, 
where $\vG(\mu) = (2n+1) \in \SP_{2n+1}$, 
\item  $\ve_{S_n} = V_{\mu}$ with $\mu = (1^n)  \lra \IC(\ol\SO_{\vG(\mu)}, \Ql)$, 
where $\vG(\mu) = (3, 2^{n-1}) \in \SP_{2n+1}$.  
\end{itemize}
\par\medskip\noindent
If $\la = (2n+1)$, $\SO_{\la}$ is the open dense orbit in $G^{\io\th}\uni$, 
hence it is similar to $\SO_1$ in Remarks 5.2 (ii).  However 
if $\la = (3,2^{n-1})$, the orbit $\SO_{\la}$ is not like $\SO_0$ in 
[loc. cit.], $\SO_{\la}$ is much bigger than the unit orbit $\SO_{\la'}$ 
with $\la' = (1^{2n+1})$. 
\par
(ii) By Corollary 12.6, $\SN_{G^{\io\th}}$ contains lots of pairs $(\SO, \SE)$
such that $\SE$ is not the constant sheaf $\Ql$. The formula 
(15.3) only involves the pairs such that $\SE = \Ql$, the constant sheaf. 
Thus the Springer correspondence is not enough to cover all the pairs
in $\SN_{G^{\io\th}}$, and we need to consider the generalized Springer correspondence.  

\par\bigskip\noindent
{\bf \S 16.  Generalized Springer correspondence for $G^{\io\th}\uni$,
             the case $H = SO_{2n+1}$ }
\addtocounter{section}{1}
\addtocounter{thm}{-3}
\par\medskip
In this section, we consider the generalized Springer correspondence for
$G^{\io\th}\uni$ with $H = SO_N$ \ ($N = 2n+1$). 

\para{16.1.}
First we will give a combinatorial description of $\SN_{G^{\io\th}}$. 
Let $\la = (\la_1, \dots, \la_N) \in \SP_N$ be a partition with $\la_N \ge 0$. 
We define a symbol $\tau = (\tau_1, \dots, \tau_N)$ of type $\la$ as follows;

\par\medskip
\begin{enumerate}
\item 
$\tau_i = \pm 1$, and $\tau_i = 1$ if $\la_i = 0$.  
\item 
$\tau_i = \tau_j$ if $\la_i = \la_j$.
\item
$\tau_k = 1$ if $\la_k$ is the largest odd number among $\la_1, \dots, \la_N$.
\end{enumerate}
\par\medskip
Then by Proposition 12.3, the set of symbols of type $\la$ is in bijection 
with the set $A_H(x)\wg$ for $x \in \SO_{\la}$.  Note that the condition (iii)
is related to the difference between $A_{G^{\th}}(x)$ and $A_H(x)$. 
Now the set $\SN_{G^{\io\th}}$ can be expressed combinatorially  
by a set $\Psi_N$ as follows; 

\begin{equation*}
\tag{16.1}
\Psi_N := \{ (\la, \tau) \mid \la \in \SP_N, \tau \text{: type $\la$} \}
        \simeq \SN_{G^{\io\th}}.
\end{equation*}
\par
An element $(\la, \tau) \in \Psi_N$ is called a {\bf cuspidal symbol }
if 
\par\medskip
\begin{enumerate}
\item \ 
$\la_i - \la_{i+1} \le 2$ for $i = 1, \dots, N$ ( here we put $\la_{N+1} = 0$ ), 
\item  \
If $\la_i - \la_{i+1} = 2$, then $\tau_i \ne \t_{i+1}$. 
\end{enumerate} 
\par\medskip
We denote by $\Psi_N^{(0)}$  the set of cuspidal symbols in $\SN_{G^{\io\th}}$.  

\para{16.2.}
Lusztig's theory of the generalized Springer correspondence for 
reductive groups  
can be extended to the case of symmetric spaces $G^{\io\th}\uni$, namely
\par\medskip
$\bullet$ The notion of a {\bf cuspidal pair}
 for $(\SO,\SE) \in \SN_{G^{\io\th}}$, 
\par\medskip
$\bullet$ The notion of { \bf induction}
 $\ind_P^G$  for a $\th$-stable parabolic subgroup $P$ and its 
$\th$-stable Levi subgroup $L$, 
\begin{equation*}
\ind_P^G : \{ \text{$L_H$-equiv. perverse sheaves on $L^{\io\th}$ } \} 
      \to \{ \text{ semisimple complexes on  $G^{\io\th}$ } \} 
\end{equation*}
The induction functor is defined by modifying the arguments in 
Section 5. 
\par
The following result was proved in [SY]. 
\begin{prop}   
Under the identification $\SN_{G^{\io\th}} \simeq \Psi_N$,
the pair $(\SO, \SE)$ is cuspidal if and only if 
the corresponding symbol $(\la, \tau) \in \Psi_N$ is cuspidal. 
\end{prop}

\remark{16.4.}
In general, $\th$-stable Levi subgroup $L$ has the form

\begin{equation*}
L_H \simeq GL_{a_1} \times G_{a_2} \times \cdots \times GL_{a_r} \times SO_{N_0},
\end{equation*}
where $N_0 = N - 2\sum_ia_i$. 
Among them, $L^{\io\th}\uni$ has a cuspidal pair $(\SO_1, \SE_1)$ only when 
$L_H \simeq (GL_1)^a \times SO_{N_0}$ for some $0 \le a \le n$.  

\para{16.5.} \
We give here some example of cuspidal pairs. It would be interesting 
to compare this with the case of reductive groups (Remark 5.2). 
\begin{enumerate}
\item  \  For $\la = (2^a, 1^{N-2a})$, the pair $(\SO_{\la}, \SE)$
is always cuspidal, for any local system $\SE$.  
\item 
In particular, for $\SO_{\la} = \{ 1\}$ with $\la = (1^N)$, \\ the pair 
$(\SO_{\la}, \Ql)$ is cuspidal 
\item 
The number $|\Psi^{(0)}_N|$ is always large. 
For example, for $N = 3, 5, 7$, we have $|\Psi^{(0)}_N| = 3, 7, 16$.  
\end{enumerate}

In general, we can determine the generating function of $|\Psi_N^{(0)}|$. 

\begin{lem} 
We define a function $q_1(n)$ for $n \ge 0$ by 
\begin{equation*}
\prod_{i = 1}^{\infty}(1 + t^i)^2 = \sum_{n \ge 0}q_1(n)t^n
\end{equation*}
Then we have $|\Psi_{2n+1}^{(0)}| = \frac{1}{2}q_1(2n + 1)$. 
\end{lem}

\para{16.7.}
We prepare some notations for formulating the generalized 
Springer correspondence. 
We define $\ScS_{G^{\io\th}}$ as the set of triples, $(L, \SO_1, \SE_1)$, 
up to the conjugation action of $H$, 
where $L$ is a $\th$-stable Levi subgroup of a $\th$-stable parabolic 
subgroup $P$, and $(\SO_1, \SE_1)$ is a cuspidal pair in $L^{\io\th}\uni$.
Since $L_H \simeq (GL_1)^a \times SO_{N_0}$ with $N_0 = N - 2a$ 
(Remark 16.4), we have
\begin{equation*}
\tag{16.2}
\ScS_{G^{\io\th}} \simeq \SC_N  := \{ (N_0, \nu, \s) \mid N \ge N_0 : \text{ odd, } 
                               (\nu,\s) \in \Psi^{(0)}_{N_0} \}
\end{equation*}

Take $\xi = (N_0, \nu,\s) \in \SC_N$ corresponding to 
$(L, \SO_1, \SE_1) \in \ScS_{G^{\io\th}}$, and 
put $K_{\xi} = \IC(\ol\SO_1, \SE_1)[\dim \SO_1]$.
Since $K_{\xi}$ is an $L_H$-equivariant perverse sheaf on $L^{\io\th}$, 
one can consider the complex $\ind_P^GK_{\xi}$ on $G^{\io\th}\uni$. 
\par
For each $\xi = (N_0, \nu,\s) \in \SC_N$, define a map 
$\vG_{\xi} : \SP_a \mapsto \Psi_N$ as follows.
Put $\wt\s = (\s_1, \dots, \s_N) = (\s_1, \dots, \s_{N_0}, 1, \dots, 1)$ 
for \ $\s = (\s_1, \dots, \s_{N_0})$. 
For $\mu \in \SP_a$, put $\la = \nu + 2\mu \in \SP_N$, and define 
\begin{equation*}
\tag{16.3}
\vG_{\xi}(\mu) = (\la, \wt\s) \in \Psi_N.
\end{equation*}

The following result gives the generalized Springer correspondence 
for $G^{\io\th}\uni$, which is an analogue of Theorem 5.1.

\begin{thm}[{[SY]}]
\begin{enumerate}  
\item 
For $\xi \in \SC_N \simeq \ScS_{G^{\io\th}}$, 
$\ind_P^G K_{\xi}$ is a semisimple perverse sheaf on $G^{\io\th}\uni$, 
equipped with $S_a$-action, and is decomposed as 

\begin{equation*}
\ind_P^G K_{\xi} = \bigoplus_{\mu \in \SP_a}V_{\mu} \otimes 
         \IC(\ol\SO,\SE)[\dim \SO],
\end{equation*} 
where $\vG_{\xi}(\mu) \lra (\SO, \SE)$ under 
the identification $\Psi_N \simeq \SN_{G^{\io\th}}$
\item
For $\xi \in \ScS_{G^{\io\th}}$, let $\SN_{G^{\io\th}}^{(\xi)}$ be 
the subset of $\SN_{G^{\io\th}}$ 
corresponding to $\{ \vG_{\xi}(\mu) \mid \mu \in \SP_a \} \subset \Psi_N$.
Then we have
\begin{equation*}
\SN_{G^{\io\th}} = \coprod_{\xi \in \ScS_{G^{\io\th}}} \SN_{G^{\io\th}}^{(\xi)}.
\end{equation*} 
\item
The correspondence $(\SO,\SE) \lra \vG_{\xi}(\mu) \lra V_{\mu}$ gives a bijection

\begin{equation*}
\SN_{G^{\io\th}} \simeq \coprod_{\xi \in \ScS_{G^{\io\th}}}S_a\wg \qquad 
                 ({\text{\bf generalized Springer correspondence}})
\end{equation*}
\end{enumerate}
\end{thm}

\par\bigskip\noindent
{\bf \S 17. Exotic symmetric spaces of higher level, the case $H = Sp_N$ }
\addtocounter{section}{1}
\addtocounter{thm}{-8}
\par\medskip
In this section, we consider the Springer correspondence 
for the exotic symmetric spaces of higher level associated to symplectic groups.
As was remarked in Section 11, the crucial difficulty in this case is that 
the number of $H$-orbits is not necessarily finite if $r \ge 3$. 
\par
Assume $H = Sp_N  \simeq Sp(V)$, and 
consider $G^{\io\th}\uni \times V^{r-1}$ ($r \ge 1$), with diagonal action of $H$. 
Let $M_1 \subset \cdots \subset M_n \subset V$ be the total isotropic flag whose 
stabilizer in $H$ is equal to $B_H$. 
Consider the variety 
\begin{align*}
   \wt\SX = \{ (x,\Bv, gB_H) \in G^{\io\th}\uni \times V^{r-1} \times \SB^H 
                  \mid g\iv xg \in B^{\io\th}, g\iv \Bv \in M_n^{r-1} \} 
\end{align*}
and define a map $\pi : \wt\SX \to G^{\io\th}\uni \times V^{r-1}$ by 
$(x,\Bv, gB_H) \mapsto (x, \Bv)$.

Let $\SQ_{n,r} = \{ \Bm = (m_1, \dots, m_r) \in \BZ_{\ge 0}^r 
                            \mid \sum m_i = n\}$, 
and put $\SQ_{n,r}^0 = \{ \Bm \in \SQ_{n,r} \mid m_r = 0 \}$. 
For given $\Bm \in \SQ_{n,r}$, define $p_1, \dots, p_{r-1}$ by 
$p_k = \sum_{i = 1}^km_i$. Put, for $\Bm \in \SQ_{n,r}^0$, 
\begin{equation*}
\SX_{\Bm} = \bigcup_{g \in H}g(U^{\io\th} \times \prod_{1 \le i <r}M_{p_i}).
\end{equation*}
\par
Let $\pi_{\Bm} : \pi\iv(\SX_{\Bm}) \to \SX_{\Bm}$ be the restriction of $\pi$ 
on $\pi\iv(\SX_{\Bm})$.  Then  
$\pi_{\Bm}$ is proper, surjective, but $\pi\iv(\SX_{\Bm})$ is not
necessarily smooth nor irreducible. 
\par
Let $W_{n,r} = S_n \ltimes (\BZ/r\BZ)^n$ be the complex reflection group.
We denote by $\SP_{n,r}$ the set of $r$-partitions 
$\Bla = (\la^{(1)}, \dots, \la^{(r)})$ such that 
$\sum_i|\la^{(i)}| = n$. 
Then it is known that 
\begin{equation*}
\tag{17.1}
W_{n,r}\wg \simeq \SP_{n,r}, \quad \wt V_{\Bla} \lra \Bla.
\end{equation*}

\par
For $\Bm \in \SQ_{n,r}^0$, we define a subset $\SP(\Bm)^0 \subset \SP_{n,r}$ by  
\begin{equation*}
\SP(\Bm)^0 = \{ \Bla \in \SP_{n,r} \mid |\la^{(i)}| = m_i \ (i \le r-2),  
       |\la^{(r-1)}| + |\la^{(r)}| = m_{r-1} \}.
\end{equation*}

The following result holds for any $r \ge 1$.  
Note that if $r = 2$, this coincides with the formula in 
Theorem  13.3.
Also note that since $\SX_{\Bm}$ has infinitely many $H$-orbits, 
we need to construct a suitable variety $X_{\Bla}$ instead of an 
$H$-orbit $\SO_{\Bla}$. 
\begin{thm} [{[S3]}]
Assume $\Bm \in \SQ^0_{n,r}$, and put $d_{\Bm} = \dim \SX_{\Bm}$.
Then $(\pi_{\Bm})_*\Ql[d_{\Bm}]$ is a semisimple perverse sheaf on $\SX_{\Bm}$, 
equipped with $W_{n,r}$-action, and is decomposed as
\begin{equation*}
(\pi_{\Bm})_*\Ql[d_{\Bm}] \simeq \bigoplus_{\Bla \in \SP(\Bm)^0}
                  \wt V_{\Bla} \otimes \IC(\ol X_{\Bla}, \Ql)[\dim X_{\Bla}],
\end{equation*}
where $X_{\Bla}$ is an $H$-stable, smooth, irreducible, locally closed 
sub-variety of $G^{\io\th}\uni \times V^{r-1}$. 
If $r \ge 3$, $X_{\Bla}$ is an infinite union of $H$-orbits. 
\end{thm}

\par
Assume $\Bm \in \SQ^0_{n,r}$, and take $\Bla \in \SP(\Bm)^0$. 
Then $X_{\Bla} \subset \SX_{\Bm}$.  Put 
\begin{equation*}
\tag{17.2}
d_{\Bla} = (\dim \SX_{\Bm} - \dim X_{\Bla})/2.
\end{equation*}
\par
For $z = (x, \Bv) \in \SX_{\Bm}$, define the {\bf Springer fibre }
$\SB^H_z \simeq \pi\iv(z)$ by 

\begin{equation*}
\tag{17.3}
\SB^H_z = \{ gB_H \in \SB^H \mid g\iv xg \in B^{\io\th}, g\iv \Bv \in M_n^{r-1} \}
\end{equation*}
\par\noindent
Then $H^i(\SB^H_z, \Ql)$ has a structure of $W_{n,r}$-module, which 
is called the Springer representation of $W_{n,r}$. 
As a corollary of Theorem 17.1, we obtain the following result, 
which is the Springer correspondence for the complex reflection group
$W_{n,r}$. 

\begin{prop} 
There exists an open dense subset $X^0_{\Bla}$ of $X_{\Bla}$ 
satisfying the following.
For $z \in X_{\Bla}^0$, 
$\dim \SB^H_z = d_{\Bla}$, and 
   {$H^{2d_{\Bla}}(\SB^H_z, \Ql) \simeq \wt V_{\Bla}$}. 
In particular, the map $X_{\Bla} \mapsto H^{2d_{\Bla}}(\SB^H_z, \Ql)$ gives a 
bijective correspondence

\begin{equation*}
\{ X_{\Bla} \mid \Bla \in \SP_{n,r} \} = 
\coprod_{\Bm \in \SQ^0_{n,r}} \{ X_{\Bla} \mid \Bla \in \SP(\Bm)^0 \} 
     \simeq \coprod_{\Bm \in \SQ^0_{n,r}} \{ \wt V_{\Bla} \mid \Bla \in \SP(\Bm)^0 \}
        = W_{n,r}\wg.
\end{equation*}
\end{prop}

\par\bigskip\noindent
{\bf \S 18.  Exotic symmetric spaces, the case $H = SO_{2n+1}$ }
\addtocounter{section}{1}
\addtocounter{thm}{-2}
\par\medskip

In this section, assume that $H = SO_{2n+1}$. Consider the variety  
$G^{\io\th}\uni \times V^{r-1}$ and the isotropic flag $(M_i)_i$ 
similarly as in Section 17.  
We fix $\xi = (N_0, \nu, \s) \lra (L, \SO_1, \SE_1) \in \SC_{G^{\io\th}}$, 
hence $L_H \simeq (GL_1)^a \times SO_{N_0}$ with $N_0 = N - 2a$.
Consider the diagram
\begin{equation*}
\tag{18.1}
\begin{CD}
L^{\io\th} @<\a<<  \wh \SX  @>\psi>>  \wt\SX @>\pi>> G^{\io\th} \times V^{r-1},
\end{CD}
\end{equation*}
where 
\begin{align*}
\wh \SX &= \{ (x, \Bv, g) \in G^{\io\th} \times V^{r-1} \times H 
               \mid g\iv xg \in P^{\io\th}, g\iv \Bv \in M_a^{r-1} \}, \\
\wt \SX &= \{ (x, \Bv, gP_H) \in G^{\io\th} \times V^{r-1} \times H/P_H
               \mid g\iv xg \in P^{\io\th}, g\iv \Bv \in M_a^{r-1} \}, 
\end{align*}
and the maps are defined as 
\begin{align*}
\psi: (x, \Bv, g) \mapsto (x, \Bv, gP_H), \quad
\pi : (x, gP_H) \mapsto x, \quad
\a : (x,\Bv, g) \mapsto \eta_P(g\iv xg),  
\end{align*}
where $\eta_P : P^{\io\th} \to L^{\io\th}$ is the natural map 
induced from the projection $P \to L \simeq P/U_P$. Note that unless 
$P = B$, the map $\wt \SX \to L^{\io\th}$ can not be defined directly.

For $\Bm \in \SQ_{a,r}$, put 

\begin{equation*}
\SX_{\Bm} = \bigcup_{g \in H}
    \biggl(\eta_P\iv(L^{\io\th}\uni) \times \prod_{1 \le i < r}M_{p_i}\biggr) 
\end{equation*}

As the restriction of the diagram (18.1) on $\SX_{\Bm}$, we have a diagram 

\begin{equation*}
\tag{18.2}
\begin{CD}
L^{\io\th} @<\a_{\Bm}<< \wh\SX_{\Bm} @>\psi_{\Bm}>> \wt\SX_{\Bm} @>\pi_{\Bm}>> \SX_{\Bm},
\end{CD}
\end{equation*}
where $\wt\SX_{\Bm} = \pi\iv(\SX_{\Bm})$, \ $\wh\SX_{\Bm} = \psi\iv(\wt\SX_{\Bm})$, 
and $\a_{\Bm}, \psi_{\Bm}, \pi_{\Bm}$ are defined as restrictions of 
$\a, \psi, \pi$, respectively. 
\par
Put $K_{\xi} = \IC(\ol\SO_1, \SE_1)[\dim \SO_1]$.  
We define a complex $\wt K_{\xi}^{(\Bm)}$ on $\wt\SX_{\Bm}$
by the condition 
\begin{equation*}
\tag{18.3}
\a_{\Bm}^*K_{\xi}[a] 
     \simeq \psi_{\Bm}^*\wt K_{\xi}^{(\Bm)}[b],
\end{equation*}
where $a$ (resp. $b$) is the dimension of the fibre of $\a_{\Bm}$
(resp. $\psi_{\Bm}$).   
We consider the complex $(\pi_{\Bm})_*\wt K^{(\Bm)}_{\xi}$ on $\SX_{\Bm}$. 
\par
Put
$
\SP(\Bm) = \{ \Bla \in \SP_{a,r} \mid |\la^{(i)}| = m_i 
                  \ (1 \le i \le r) \}.
$
Since the number of $H$-orbits is infinite for $G^{\io\th}\uni$, we need to
construct subvarieties $X_{\Bla}$ as in Section 17.

\begin{prop}  
Let $\Bm \in \SQ_{a,r}$.  For each $\Bla \in \SP(\Bm)$, 
one can define a subvariety $X_{\Bla}$ of $\SX_{\Bm}$ and 
an $H$-orbit $\SO_{[\Bla]}$ in $G^{\io\th}\uni$ satisfying 
the following; 
\par\smallskip
\begin{enumerate}
\item 
$X_{\Bla}$ is a locally closed, smooth, irreducible, $H$-stable
subset of $\SX_{\Bm}$.
\item
The map $(x,\Bv) \mapsto x$ gives a locally trivial fibration 
$f_{\Bla} : X_{\Bla} \to \SO_{[\Bla]}$, \\
whose fibre is isomorphic to 
an open dense subset of $\prod_{1 \le i < r}M_{p_i}$.  
\item
If $r = 1$, $X_{\Bla}$ coincides with $\SO_{[\Bla]} = \SO_{2\la^{(1)} + \nu}$. 
\end{enumerate}
\end{prop}

For $\SO = \SO_{[\Bla]}$, there exists a unique pair 
$(\SO,\SE) \in \SN^{(\xi)}_G$ appearing in $\ind_P^GK_{\xi}$. 
We define a local sytem $\SE_{\Bla}$ on $X_{\Bla}$ by 
$\SE_{\Bla} = f_{\Bla}^*\SE$. 

\begin{thm}[{[S4]}]  
For each $\Bm \in \SQ_{a,r}$, 
$(\pi_{\Bm})_*\wt K_{\xi}^{(\Bm)}$ is a semisimple perverse sheaf 
on $\SX_{\Bm}$, equipped with $W_{a,r}$-action, and is decomposed as
\begin{equation*}
(\pi_{\Bm})_*\wt K_{\xi}^{(\Bm)} \simeq \bigoplus_{\Bla \in \SP(\Bm)}
               \wt V_{\Bla}\otimes \IC(\ol X_{\Bla}, \SE_{\Bla})[\dim X_{\Bla}].
\end{equation*}
In particular, for each $\xi \in \ScS_{G^{\io\th}}$, we have a bijection 
\begin{equation*}
\{ X_{\Bla} \mid \Bla \in \SP_{a,r} \} 
   = \coprod_{\Bm \in \SQ_{a,r}}\{ X_{\Bla} \mid \Bla \in \SP(\Bm) \}
        \simeq \coprod_{\Bm \in \SQ_{a,r}}\{ \wt V_{\Bla} \mid \Bla \in \SP(\Bm) \}
        = W_{a,r}\wg 
\end{equation*} 
\end{thm}

\remark{18.4.} \
The construction of $(\pi_{\Bm})_*\wt K^{(\Bm)}_{\xi}$ is essentially 
the same as the construction of the induction $\ind_P^GK_{\xi}$ in 
Section 16, which is a variant of the induction functor of Lusztig. 
But in contrast to the cases of reductive groups or of $G^{\io\th}\uni$ 
for $H = SO_{2n+1}$, $(\pi_{\Bm})_*\wt K^{(\Bm)}_{\xi}$ is not a priori 
a semisimple complex.  So the previous discussion can not be applied 
directly, and we need a special care.

\par\bigskip\noindent
{\bf \S 19. Symmetric spaces in characteristic 2 }
\addtocounter{section}{1}
\addtocounter{thm}{-3}
\par\medskip
The discussion in Section 11 makes sense even if $p = 2$, and 
one can define a symmetric space $G^{\io\th}\uni$ associated to
classical groups in characteristic 2. In this case, a different phenomenon 
appears since the fundamental properties given in Section 10 do not
hold. However, an analogue of the Springer correspondence still holds for
them. The Springer correspondence for $G^{\io\th}\uni$ was studied in 
[DSY], which we briefly explain below. 
\par 
Let $V'$ be an $N$-dimensional vector space over $\Bk$ with $p = 2$, 
and consider $G = GL(V')$. The involution $\th : G \to G$ is defined as in 
Section 11, and we consider $G^{\th}$ and $G^{\io\th}$. 
Here $G^{\th} \simeq Sp(V)$, where 
$V$ is an $2n$-dimensional subspace of $V'$ if $N = 2n+1$ and $V = V'$ 
if $N = 2n$. Put $H = Sp(V)$ and $\Fh = \Lie H$.  Then 
$\Fh = \Fs\Fp(V)$: the symplectic Lie algebra with $p = 2$.  
$\th$ induces an involution $\th$ on $\Fg = \Lie G$, and we consider 
the subalgebra  
$\Fg^{\th} = \{ x \in \Fg \mid \th(x) = x \}$. 
The case where $N$ is even, the situation is rather simple, namely 
we have

\begin{prop} 
Assume that $N = 2n$.  Then 
\begin{enumerate} 
\item 
$\Fg^{\th} \simeq \Lie H = \Fs\Fp(V)$: the symplectic Lie algebra over $\Bk$.
\item
$G^{\io\th}\uni \simeq \Fg^{\th}\nil$, compatible with the action of $H$.
\end{enumerate}
\end{prop} 
By the above result, considering the Springer correspondence for $G^{\io\th}\uni$ 
is equivalent to considering it for $\Fs\Fp(V)\nil$. 
The Springer correspondence for the nilpotent cone of the Lie algebras 
in characteristic 2 was established by T. Xue [X].  
Thus the Springer 
correspondence for $G^{\io\th}\uni$ follows from her result.  
In particular, $G^{\io\th}\uni/\!\sim_H$ is in bijection with $\SP_{n,2}$, 
via $\SO_{\Bla} \lra \Bla$, and the following formula holds.

\begin{equation*}
\tag{19.1}
\pi_*\Ql[\dim G^{\io\th}\uni] \simeq 
   \bigoplus_{\Bla \in \SP_{n,2}}\wt V_{\Bla}\otimes \IC(\ol\SO_{\Bla}, \Ql)[\dim \SO_{\Bla}],
\end{equation*} 
where $\wt V_{\Bla}$ is the irreducible representation of $W_{n,2}$ corresponding 
to $\Bla \in \SP_{n,2}$. 
\par
Next consider the case where $N$ is odd. In this case, a quite different 
situation occurs. 
We have

\begin{prop} 
Assume that $N = 2n+1$.  Then 
\begin{enumerate}
\item
$G^{\io\th}\uni \simeq \Fg^{\th}\nil$, compatible with the action of $H \simeq G^{\th}$. 
\item 
$\Fg^{\th}\simeq \Fh \times V \times \Bk$, where $H$ acts on $\Fh \times V$ diagonally, 
and acts on $\Bk$ trivially. 
\end{enumerate}
\end{prop}

\remarks{19.3.} 
(i) \ $\Lie G^{\th} \ne \Fg^{\th}$. \\
(ii) \ The number of $H$-orbits in $\Fg^{\th}\nil$ is infinite, hence 
the number of $G^{\th}$-orbits in $G^{\io\th}\uni$ is infinite. 
(This gives a counter example for Proposition 10.2 in the case $p = 2$.)  \\
(iii) \ The situation is similar to the exotic case, namely 
the action of $H$ on $G^{\io\th}$ corresponds to 
the diagonal action of $H$ on $\Fh \times V$.  
But the nilpotent part $\Fg^{\th}\nil$ is not exactly the same as 
$\Fh\nil \times V$.  In our case, we have 
$\Fg^{\th}\nil = (\Fh\nil \times V) \cap \Fg\nil$ for $\Fg = \Lie G = \Fg\Fl(V')$.   
\par\medskip
The Springer correspondence for $G^{\io\th}$ was established in [DSY].
Note that in this case, it occurs a very similar phenomenon as in the case of 
exotic symmetric space of level $r = 3$.  In fact, since $G^{\io\th}\uni$ has 
infinitely many $H$-orbits, we need to consider the varieties such as $X_{\Bla}$ 
in Section 17.  In this case, $X_{\Bla}$ is parametrized by $\Bla \in \SP_{n,3}$, 
and the Springer correspondence is described by using the complex reflection group
$W_{n,3}$. 
\par
The reason why $W_{n,3}\wg$ appears in the Springer correspondence is 
explained as follows.  It is proved by Kato [K2], by using the deformation 
arguments on the schemes over $\BZ$,  that the Springer correspondence 
for the exotic symmetric space 
$\Fg^{-\th}\nil \times \BV$ for $H = Sp(\BV)$ in the case $p \ne 2$ 
can be interpreted by  
the Springer correspondence for $\Fs\Fp(V)\nil$ in the case $p = 2$.
Thus, as a generalization of his result, one can expect that 
\begin{equation*}
\tag{19.2}
\Fs\Fp(V)\nil \times V  \approx \Fg^{-\th}\nil \times \BV^2 \approx 
                           G^{\io\th}\uni \times \BV^2.  
\end{equation*}
In [DSY] it is proved that this certainly holds. Note that 
the right hand side is the exotic symmetric space of $r = 3$. 

\par

\par\bigskip\bigskip
\noindent
T. Shoji \\
School of Mathematical Sciences, Tongji University \\ 
1239 Siping Road, Shanghai 200092, P. R. China  \\
E-mail: \verb|shoji@tongji.edu.cn|

\end{document}